\documentclass[a4paper,11pt]{amsart}
\usepackage{fancyhdr}
\usepackage{amsmath}
\usepackage{dsfont}
\usepackage{hyperref}
\usepackage[mathscr]{eucal}
\usepackage[cp1251]{inputenc}
\usepackage[english]{babel}
\usepackage{cite,enumerate,float,indentfirst}
\usepackage{graphicx}
\usepackage{xcolor}
\usepackage{latexsym,a4,mathrsfs,amsthm,amsmath,amssymb,url}

\numberwithin{equation}{section}
\newtheorem{theorem}{Theorem}[section]
\newtheorem{lemma}[theorem]{Lemma}
\newtheorem{statement}[theorem]{Statement}
\newtheorem{definition}[theorem]{Definition}

\newtheorem{remark}[theorem]{Remark}
\newcommand{\R}{\mathbb{R}}
\newcommand{\C}{\mathbb{C}}
\newcommand{\T}{\mathbb{T}}
\newcommand{\Levy}{{\EuScript L}}
\newcommand{\sgn}{\text{sgn}}

\newcommand{\Cond}{{\bf (C0)}}

\newcommand{\X}{{\bf X}}
\newcommand{\A}{{\bf A}}
\newcommand{\M}{{\bf M}}
\newcommand{\B}{{\bf B}}
\newcommand{\Q}{{\bf Q}}
\newcommand{\EE}{{\bf E}}
\newcommand{\F}{{\bf F}}
\newcommand{\G}{{\bf G}}
\newcommand{\HH}{{\bf H}}
\newcommand{\OO}{{\bf O}}
\newcommand{\OP}{{\bf P}}
\newcommand{\I}{{\bf I}}
\newcommand{\RR}{{\bf R}}
\newcommand{\V}{{\bf V}}
\newcommand{\J}{{\bf J}}

\DeclareMathOperator{\Tr}{Tr}
\DeclareMathOperator{\Rank}{rk}
\DeclareMathOperator{\dist}{dist}
\DeclareMathOperator{\Sp}{span}
\DeclareMathOperator{\E}{\mathbb{E}}
\DeclareMathOperator{\Pb}{\mathbb{P}}
\DeclareMathOperator{\one}{\mathds{1}}

\DeclareMathOperator{\imag}{Im}
\DeclareMathOperator{\spr}{spread}
\DeclareMathOperator{\Real}{Re}

\begin{document}

\vspace{1in}

\title
{\bf ELLIPTIC LAW FOR REAL RANDOM MATRICES}

%\vspace{2in}

\author[A. Naumov]{Alexey Naumov$^1$}
\address{A. Naumov\\
 Faculty of Computational Mathematics and Cybernetics\\
 Moscow State University \\
 Moscow, Russia
 }
\email{naumovne@gmail.com, anaumov@math.uni-bielefeld.de}
\thanks{$^1$This research was supported by DAAD}

%\subjclass[2000]{60F05, 15A52, 20B30, 60B15, 60C05, 60F05, 60G50, 60E10, 20C30, 32F18}

\date{\today}

\begin{abstract}
In this paper we consider ensemble of random matrices $\X_n$ with independent identically distributed vectors $(X_{ij}, X_{ji})_{i \neq j}$ of entries. Under assumption
of finite fourth moment of matrix entries it is proved that empirical spectral distribution of eigenvalues converges in probability
to a uniform distribution on the ellipse. The axis of the ellipse are determined by correlation between $X_{12}$ and $X_{21}$. This result is called Elliptic Law.
Limit distribution doesn't depend on distribution of matrix elements and the result in this sence is universal.
\end{abstract}

\maketitle
\tableofcontents

%\selectlanguage{english}

\newpage

\section{Introduction}

Let us consider real random matrix $\X_n(\omega) = \{X_{i j}(\omega)\}_{i, j=1}^n$ and assume that the following conditions $\Cond$ hold\\
a) Pairs $(X_{ij}, X_{ji}), i \neq j$ are independent identically distributed (i.i.d.) random vectors;\\
b) $\E X_{1 2} = \E X_{2 1} = 0, \E X_{1 2}^2 = \E X_{2 1}^2 = 1$ and $\max(\E|X_{1 2}|^4, \E|X_{2 1}|^4) \le M_4$;\\
c) $\E ( X_{1 2} X_{2 1} ) = \rho$, $|\rho| \le 1$;\\
d) The diagonal entries $X_{ii}$ are i.i.d. random variables, independent of off-diagonal entries, $\E X_{1 1} =0$ and $\E X_{1 1}^2 < \infty$.

Denote by $\lambda_1, ..., \lambda_n$ the eigenvalues of the matrix $n^{-1/2} \X_n$ and define empirical spectral measure by
$$
\mu_n(B) = \frac{1}{n} \# \{1 \le i \le n: \lambda_i \in B \}, \quad B \in \mathcal{B}(\C),
$$
where $\mathcal{B}(\C)$ is a Borel $\sigma$-algebra of $\C$.

We say that the sequence of random probability measures $m_n(\cdot)$
converges weakly in probability to probability measure $m(\cdot)$ if for all continues and bounded functions $f: \C \rightarrow \C$ and all $\varepsilon > 0$
$$
\lim_{n \rightarrow \infty}\Pb \left ( \left | \int_\C f(x) m_{n}(dz) - \int_\C f(x)m(dz) \right | > \varepsilon \right ) = 0.
$$

We denote weak convergence by symbol $\xrightarrow{weak}$.

A fundamental problem in the theory of random
matrices is to determine the limiting distribution of $\mu_n$
as the size of the random matrix tends to infinity. The main result of this paper is the following
\begin{theorem} {\bf (Elliptic Law)}\label{th:main}
Let $\X_n$ satisfies condition $\Cond$ and $|\rho| < 1$.   Then $\mu_n \xrightarrow{weak} \mu$ in probability, and $\mu$ has a density $g$:
$$
g(x, y) = \begin{cases}
  \frac{1}{\pi (1 - \rho^2)}, & x, y \in \EuScript E, \\
  0, & \text{otherwise,}
\end{cases}
$$
where
$$
\EuScript E := \left \{ x,y \in \R: \frac{x^2}{(1+\rho)^2} + \frac{y^2}{(1-\rho)^2} \le 1 \right \}.
$$
\end{theorem}

Theorem~\ref{th:main} asserts that under assumption of finite fourth moment empirical distribution weakly converges in probability to uniform distribution on the ellipse.
The axis of the ellipse are determined by correlation between $X_{12}$ and $X_{21}$. This result was called by Girko ``Elliptic Law''.
Limit distribution doesn't depend on distribution of matrix elements and the result in this sense is universal.

\begin{figure}
\begin{center}
\scalebox{.4}{\includegraphics{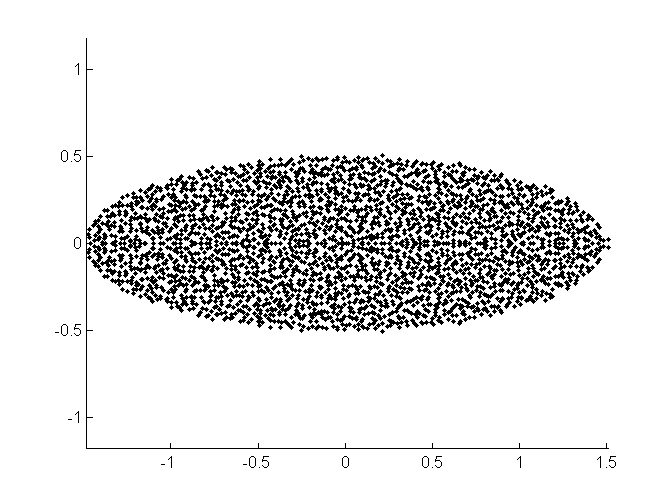}}
\scalebox{.4}{\includegraphics{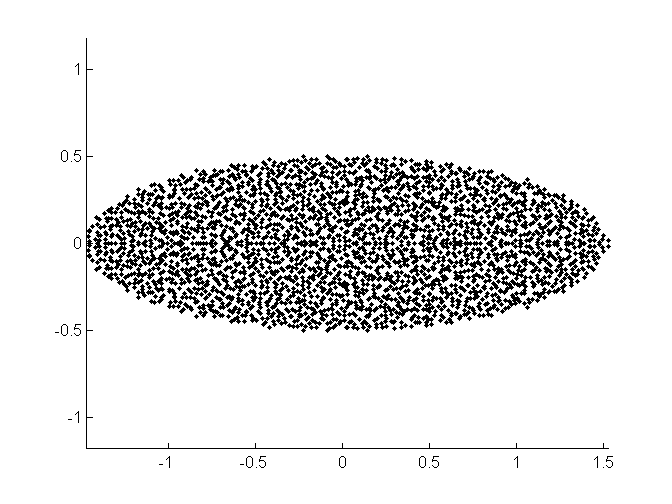}}
\end{center}
\caption{Eigenvalues of matrix $n^{-1/2} \X$ for $n = 3000$ and $\rho = 0.5$.
On the left, each entry is an iid Gaussian normal random variable.
On the right, each entry is an iid Bernoulli random variable, taking the values $+1$ and $-1$ each with probability $1/2$. }
\label{fig:ellip1}
\end{figure}

\begin{figure}
\begin{center}
\scalebox{.4}{\includegraphics{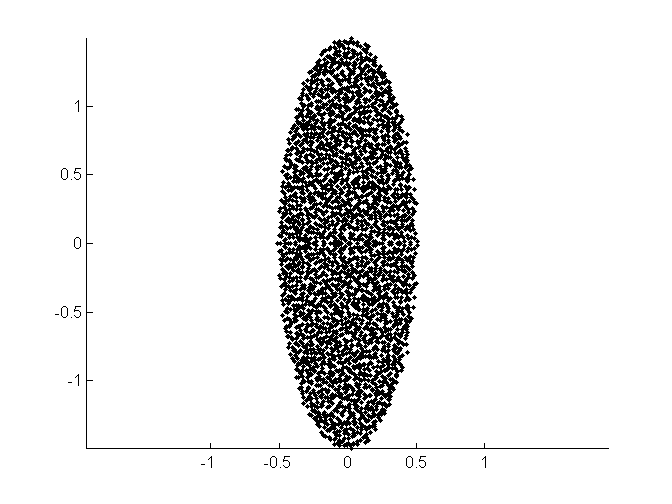}}
\scalebox{.4}{\includegraphics{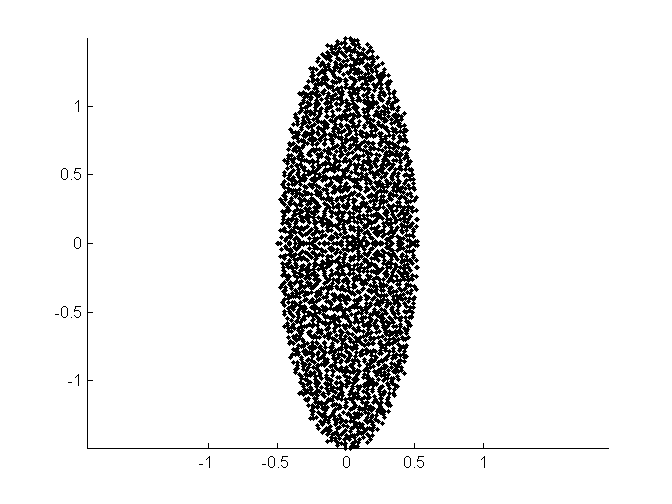}}
\end{center}
\caption{Eigenvalues of matrix $n^{-1/2} \X$ for $n = 3000$ and $\rho = -0.5$.
On the left, each entry is an iid Gaussian normal random variable.
On the right, each entry is an iid Bernoulli random variable, taking the values $+1$ and $-1$ each with probability $1/2$. }
\label{fig:ellip2}
\end{figure}

Figure~\ref{fig:ellip1} illustates Elliptic law for $\rho = 0.5$ and Figure~\ref{fig:ellip2} -- for $\rho = -0.5$.

%Figure~\ref{fig:ellip and circ} illustates Elliptic law.

In 1985 Girko proved elliptic law for rather general ensembles of random matrices under assumption
that matrix elements have a density, see~\cite{Girko1985} and~\cite{Girko2006}.
Girko used method of characteristic functions. Using $V$-transform he reduced problem
to the problem for Hermitian matrices $(n^{-1/2} \X_n - z \I)^{*}(n^{-1/2} \X_n - z \I)$ and
established convergence of empirical spectral distribution of singular values of $n^{-1/2}\X_n - z \I$ to the
limit which determines the elliptic law.

Let elements of real asymmetric random matrix $\X$ have Gaussian distribution with zero mean and correlations
$$
\E X_{i j}^2 = 1 \text{ and } \E X_{i j} X_{i j} = \rho, \quad i \neq j,  \quad |\rho| < 1.
$$
The ensemble of such matrices can be specified by the probability measure
$$
\Pb(dX) \thicksim \exp \left [-\frac{n}{2(1 - \rho^2)} \Tr(X X^T - \rho X^2) \right ].
$$
It was proved that $\mu_n \xrightarrow{weak} \mu$, where $\mu$ has a density from Theorem~\ref{th:main}, see~\cite{SomStein1988}. We will use this result to prove Theorem~\ref{th:main} in the general case.
\begin{remark}
This result can be generalized to an ensemble of Gaussian complex asymmetric matrices. In this case, the invariant measure is
$$
\Pb(dX) \thicksim \exp \left [-\frac{n}{1 - |\rho|^2} \Tr(X X^T - 2 \Real \rho X^2) \right ]
$$
and $\E|X_{i j}|^2 =1, \E X_{i j} X_{j i} = |\rho| e^{2 i \theta}$ for $i \neq j$. Then the limit measure has a uniform density inside an ellipse which is centered at zero and has semiaxes $1 + |\rho|$ in the direction $\theta$ and $1 - |\rho|$ in the direction $\theta + \pi/2$.
\end{remark}
For the discussion of elliptic law in Gaussian case see also~\cite{FyodSomm1998},~\cite[Chapter~18]{Akem2011} and~\cite{Ledoux2008}.

We repeat physical motivation of models of random matrices which satisfy condition $\Cond$ from~\cite{SomStein1988}:
{\it "The statistical properties of random asymmetric matrices may be important in the understanding of the behavior of
certain dynamical systems far from equilibrium. One example is the dynamics of neural networks.
A simple dynamic model of neural network consists of $n$ continues "scalar" degrees of freedom("neurons")
obeying coupled nonlinear differential equations ("circuit equations"). The coupling between the neurons
is given by a synaptic matrix $\X$ which, in general, is asymmetric and has a substantial degree of disorder.
In this case, the eigenstates of the synaptic matrix play an important role in the dynamics particulary when the neuron nonlinearity is not big"}.

It will be interesting to prove Theorem~\ref{th:main} only under assumption of finite second moment and prove sparse analogs. It is the direction of our further research.

If $\rho = 0$ we assume that all entries of $\X_n$ are independent random variables and Circular law holds (see~\cite{BaiSilv2010},\cite{TaoVu2010},\cite{GotTikh2010}):
\begin{theorem} {\bf (Circular law)}
Let $\X_n$ be a random matrix with independent identically distributed entries, $\E X_{i j} =0$ and $\E X_{i j}^2 = 1$. Then $\mu_n \xrightarrow{weak} \mu$ in probability, and $\mu$ has uniform density on the unit circular.
\end{theorem}

\begin{figure}
\begin{center}
\scalebox{.4}{\includegraphics{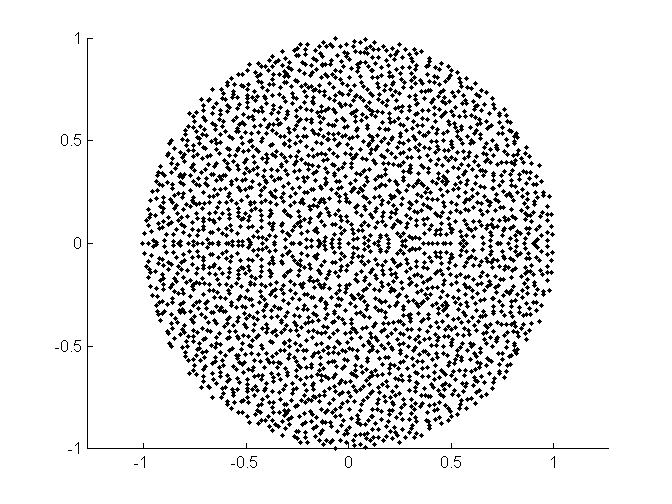}}
\scalebox{.4}{\includegraphics{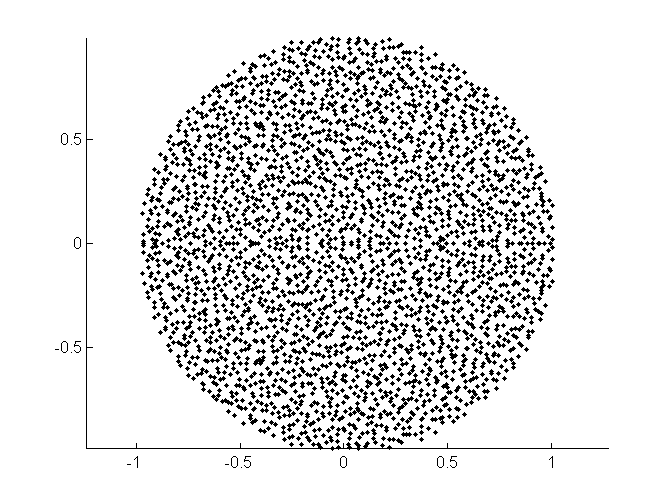}}
\end{center}
\caption{Eigenvalues of matrix $n^{-1/2} \X$ for $n = 3000$ and $\rho = 0$.
On the left, each entry is an iid Gaussian normal random variable.
On the right, each entry is an iid Bernoulli random variable, taking the values $+1$ and $-1$ each with probability $1/2$. }
\label{fig:circ}
\end{figure}

See Figure~\ref{fig:circ} for illustration of Circulaw law.

If $\rho = 1$ then matrix $\X_n$ is symmetric and its eigenvalues are real numbers. In this case the next theorem is known as a Wigner's semi-circular law (see~\cite{BaiSilv2010}):
\begin{theorem} {\bf (Semi-circular law)}
Let $\X_n$ be a symmetric random matrix with independent identically distributed entries for $i \geq j$, $\E X_{i j} =0$, $\E X_{i j}^2 = 1$. Then $\mu_n \xrightarrow{weak} \mu$ in probability, and $\mu$ has a density $g$:
$$
g(x) = \begin{cases}
  \frac{1}{2 \pi} \sqrt{4 - x^2}, & -2 \le x \le 2, \\
  0, & \text{ otherwise.}
\end{cases}
$$
\end{theorem}

Throughout this paper we assume that all random variables are defined on common probability space $(\Omega, {\EuScript F}, \Pb)$
and we will write almost surely (a.s) instead of $\Pb$-almost surely. By $\Tr(\A)$ and $\Rank(\A)$ we mean trace and rank of the matrix $\A$ respectively. We denote
singular values of matrix $\A$ by $s_i(\A)$ and $s_1(\A) \geq s_2(\A) \geq ... \geq s_n(\A)$.
For vector $x = (x_1, ... , x_n)$ we introduce $||x||_2 := (\sum_{i=1}^n x_i^2 )^{1/2}$ and $||x||_3 := (\sum_{i=1}^3 |x_i|^3)^{1/3}$.
We denote unit sphere and unit ball by $S^{n-1} := \{x: ||x||_2 = 1 \}$ and $B_1^n := \{x: ||x||_2 \le 1\}$ respectively.
For matrix $\A$ define spectral norm by $||\A|| := \sup_{x: ||x||_2 = 1} ||\A x||_2$ and Hilbert-Schmidt norm by $||\A||_{HS} := (\Tr(\A^*\A))^{1/2}$.
By $[n]$ we mean the set $\{1, ..., n\}$. We denote by $\mathcal{B}(\mathbb T)$ - Borel $\sigma$-algebra of $\mathbb T$, where $\mathbb T = \R$ or $\C$.

\section{Proof of the main result}
Further we will need the definition of logarithmic potential (see~\cite{SaffTotik1997}) and uniform integrability of function with respect to the sequence of probability measures.
\begin{definition}
The logarithmic potential $U_m$ of measure $m(\cdot)$ is a function $U_m: \C \rightarrow (-\infty, + \infty]$ defined for all $z \in \C$ by
$$
U_m (z) = -\int_\C \log|z - w| m(dw).
$$
\end{definition}
\begin{definition}
The function $f: \T \rightarrow \R$, where $\T = \C$ or $\T = \R$, is uniformly integrable in probability with respect to the sequence of
random measures $\{m_n\}_{n \geq 1}$ on $(\T, \mathcal B(\T))$ if for all $\varepsilon>0$:
$$
\lim_{t \rightarrow \infty} \varlimsup_{n \rightarrow \infty} \mathbb P \left (\int_{|f|>t} |f(x)| m_n(dx) > \varepsilon \right ) = 0.
$$
\end{definition}

Let $s_1(n^{-1/2}\X - z \I) \geq s_2(n^{-1/2}\X - z \I) \geq ... \geq s_n(n^{-1/2}\X - z \I)$ be singular values of $n^{-1/2} \X_n - z \I$ and
$$
\nu_n(z, B) = \frac{1}{n} \# \{i \geq 1: s_i(n^{-1/2}\X - z\I) \in B \}, \quad B \in \mathcal{B}(\R) -
$$
empirical spectral measure of singular values. We will omit argument $z$ in notation of measure $\nu_n(z,B)$ if it doesn't confuse.

The convergence in the Theorem~\ref{th:main} will be proved via convergence of logarithmic potential of $\mu_n$ to the logarithmic potential of $\mu$.
We can rewrite logarithmic potential of $\mu_n$ via the logarithmic moments of measure $\nu_{n}$ by
\begin{align*} \label{eq:log. potential}
&U_{\mu_n} (z) = -\int_\C \log|z - w| \mu_n(dw) = -\frac{1}{n} \log \left |\det \left (\frac{1}{\sqrt n} \X_n - z\I \right )\right| \\
&= - \frac{1}{2 n} \log \det \left (\frac{1}{\sqrt n} \X_n-z\I \right )^{*}\left (\frac{1}{\sqrt n} \X_n-z\I \right ) =
-\int_0^\infty \log x \nu_{n}(dx).
\end{align*}
This allows us to consider Hermitian matrix $(n^{-1/2} \X_n-z\I)^{*}(n^{-1/2} \X_n-z\I)$ instead of asymmetric $n^{-1/2} \X$. To prove Theorem~\ref{th:main} we need the following
\begin{lemma} \label{l:Girko}
Let $(\X_n)_{n \geq 1}$ be a sequence of $n \times n$ random matrices. Suppose that for a.a. $z \in \mathbb C$ there exists a probability measure $\nu_z$ on $[0, \infty)$ such that\\
a) $\nu_n \xrightarrow{weak} \nu_z$ as $n \rightarrow \infty$ in probability\\
b) $\log$ is uniformly integrable in probability with respect to $\{\nu_n\}_{n \geq 1}$.

Then there exists a probability measure $\mu$ such that\\
a) $\mu_{n} \xrightarrow{weak} \mu$ as $n \rightarrow \infty$ in probability \\
b) for a.a. $z \in \C$
$$
U_{\mu}(z) = - \int_{0}^\infty \log x \nu_z(dx).
$$
\end{lemma}
\begin{proof}
See \cite[Lemma~4.3]{BordCh2011} for the proof.
\end{proof}

\begin{proof} (Proof of Theorem~\ref{th:main})
Our aim is to prove convergence of $\nu_n$ to $\nu_z$, uniform integrability of $\log( \cdot )$ with respect to $\{\nu_n\}_{n \geq 1}$ and show that $\nu_z$ determines elliptic law.

From Theorem~\ref{th: log uniform integr} we can conclude uniform integrability of $\log(\cdot)$. The proof of Theorem~\ref{th: log uniform integr} is based on Theorem~\ref{th:smallest singular value} and
some additional results.

In Theorem~\ref{th:Stieltjes transform} it is proved that $\nu_n \xrightarrow{weak} \nu_z$ in probability, where $\nu_z$ is some probability measure, which doesn't depend on distribution of elements of matrix $\X$.

If matrix $\X$ has Gaussian elements we redenote $\mu_n$ by $\hat \mu_n$.

By Lemma~\ref{l:Girko} there exists probability measure $\hat \mu$ such that $\mu_n \xrightarrow{weak} \hat \mu$ in probability and $U_{\hat \mu}(z) = - \int_{0}^\infty \log x \nu_z(dx)$. But in Gaussian case
$\mu_n \xrightarrow{weak} \mu$ in probability and $U_{\mu}(z) = - \int_{0}^\infty \log x \nu_z(dx)$. We know that $\nu_z$ is the same for all matrices which satisfy condition $\Cond$ and we have
$$
U_{\hat \mu}(z) = - \int_{0}^\infty \log x \nu_z(dx) = U_{\mu}(z).
$$
From unicity of logarithmic potential we conclude that $\hat \mu = \mu$.
\end{proof}

\section{Least singular value} \label{sec:least singular value}

From properties of the largest and the smallest singular values
$$
s_1(\A) = ||\A|| = \sup_{x: ||x||_2 = 1} ||\A x||_2, \quad s_n(\A) = \inf_{x: ||x||_2 = 1} ||\A x||_2.
$$
To prove uniform integrability of $\log(\cdot)$ we need to estimate probability of the event $\{s_n(\A) \le \varepsilon n^{-1/2}, ||\X|| \le K \sqrt n \}$,
where $\A = \X - z \I$. We can assume that $\varepsilon n^{-1/2} \le K n^{1/2}$. If $|z| \geq 2 K \sqrt n$ then probability of the event is automatically zero.
So we can consider the case when $|z| \le 2 K n^{1/2}$. We have $||\A|| \le ||\X|| + |z| \le 3 K n^{1/2}$.
In this section we prove theorem
\begin{theorem} \label{th:smallest singular value}
Let $\A = \X - z \I$, where $\X$ is $n \times n$ random matrix satisfying $\Cond$. Let $K > 1$. Then for every $\varepsilon >0$ one has
$$
\Pb (s_n(\A) \le \varepsilon n^{-1/2}, ||\A|| \leq 3 K \sqrt n) \le C(\rho) \varepsilon^{1/8} + C_1(\rho) n^{-1/8},
$$
where $C(\rho), C_1(\rho)$ are some constants which can depend only on $\rho, K$ and $M_4$.
\end{theorem}
\begin{remark}
Mark Rudelson and Roman Vershynin in~\cite{RudVesh2008} and Roman Vershynin in~\cite{Veshyn2011} found bounds for the least singular value of matrices with independent entries and symmetric matrices respectively. In this section we will follow their ideas.
\end{remark}

\subsection{The small ball probability via central limit theorem} \label{sec:SBP CLT}

We recall definition of Levy concentration function
\begin{definition}
Levy concentration function of random variable $Z$ with values from $\R^d$ is a function
$$
\Levy (Z, \varepsilon) = \sup_{v \in \R^d}\Pb(||Z - v||_2 < \varepsilon ).
$$
\end{definition}
The next statement gives the bound for Levy concentration function of sum of independent random variables in $\R$.
\begin{statement}\label{st:SBP CLT}
Let $\{a_i \xi_i + b_i \eta_i\}_{i \geq 1}$ be independent random variables, $\E \xi_i = \E \eta_i = 0$, $\E \xi_i^2 \geq 1, \E \eta_i^2 \geq 1$,
$\E \xi_i \eta_i = \rho$, $\max(\E \xi_i^4, \E \eta_i^4) \le M_4$, $a_i^{-1}b_i = O(1)$. We assume that $\tau (2n)^{-1/2} \le |a_i| \le (\delta n)^{-1/2}$,
where $\delta, \tau$ are some constants. Then
$$
\Levy \left (\sum_{i=1}^n (a_i \xi_i + b_i \eta_i), \varepsilon \right )\le \frac{C \varepsilon}{(1- \rho^2)^{1/2}} + \frac{C_1}{(1-\rho^2)^{3/2} n^{1/2}}.
$$
\end{statement}
\begin{proof}
Set $\sigma_{i1}^2 = \E \xi_i^2$ and $\sigma_{i2}^2 = \E \eta_i^2$. It is easy to see that
\begin{align*}
&\sigma^2 = \E (\sum_{i=1}^n Z_i)^2 =\\
&=\sum_{i=1}^n |a_i|^2 (\sigma_{i1}^2 + 2 \rho \sigma_{i1} \sigma_{i2} a_i^{-1}b_i + \sigma_{i2}^2(a_i^{-1}b_i)^2)
\geq (1-\rho^2) \sum_{i=1}^n \sigma_{i1}^2 |a_i|^2
\end{align*}
and
$$
\sum_{i=1}^n \E|a_i \xi_i + b_i \eta_i|^3 \le \sum_{i=1}^n |a_i|^3 \E|\xi_i + a_i^{-1} b_i \eta_i|^3 \le C' M_4 ||a||_3^3,
$$
where we have used the fact $a_i^{-1} b_i = O(1)$. By Central Limit Theorem~\ref{th:CLT} for arbitrary vector $v \in \R$
$$
\Pb \left ( \left|\sum_{i = 1}^n (a_i \xi_i + b_i \eta_i) - v \right | \le \varepsilon \right )
\le \mathbb P \left ( \left| g' - v \right| \le \varepsilon \right ) + C'' \frac{\sum_{i=1}^n \E|a_i \xi_i + b_i \eta_i|^3}{\sigma^3},
$$
where $g'$ has gaussian distribution with zero mean and variance $\sigma^2$. The density of $g'$ is uniformly bounded by $1/\sqrt{2 \pi \sigma^2}$.
We have
$$
\mathbb P \left (|\sum_{i = 1}^n (a_i \xi_i + b_i \eta_i) - v| \le \varepsilon \right ) \le \frac{C \varepsilon}{(1- \rho^2)^{1/2}} + \frac{C_1}{(1-\rho^2)^{3/2} n^{1/2}}.
$$
We can take maximum and conclude the statement.
\end{proof}
\begin{remark}\label{r:SBP CLT}
Let us consider the case $b_i = 0$ for all $i \geq 1$. It is easy to show that
$$
\Levy \left (\sum_{i=1}^n a_i \xi_i, \varepsilon \right )\le C (\varepsilon +  n^{-1/2}).
$$
\end{remark}

\subsection{Decomposition of the sphere and invertibility} \label{sec: decomposition}

To prove Theorem~\ref{th:smallest singular value}, we shall partition the unit sphere $S^{n-1}$
into the two sets of compressible and incompressible vectors, and show the invertibility of $\A$ on each set separately.

\begin{definition}(Compressible and incompressible vectors)
Let $\delta, \tau \in (0, 1)$. A vector $x \in \mathbb R^n$ is called sparse if $|supp(x)| \le \delta n$.
A vector $x \in S^{n-1}$ is called compressible if $x$ is within Euclidian distance $\tau$ from the set of all sparse vectors.
A vector $x \in S^{n-1}$ is called incompressible if it is not compressible.
The sets of sparse, compressible and incompressible vectors will be denoted by Sparse = Sparse ($\delta$), Comp = Comp ($\delta, \tau$) and Incomp = Incomp($\delta, \tau$) respectively.
\end{definition}

%\subsection{Invertibility for the compressible vectors} \label{sec: comp}

We first estimate $||\A x||$ for a fixed vector $x \in S^{n-1}$. The next statement can be found in~\cite{Veshyn2011}

\begin{lemma}\label{l:comp. vectors}
Let $\A$ be a matrix from Theorem~\ref{th:smallest singular value} and let $K > 1$. There exist constants $\delta, \tau, c \in (0, 1)$ that depend only on $K$ and $M_4$
and such that the following holds. For every $u \in \R^n$, one has
\begin{equation}\label{eq:comp. vectors}
\Pb \left( \inf_{\frac{x}{||x||_2} \in Comp(\delta, \tau) } ||\A x - u||_2 / ||x||_2 \le c_4 \sqrt n, ||\A|| \le 3 K \sqrt n \right ) \le 2 e^{-c n}.
\end{equation}
\end{lemma}
\begin{proof}
See~\cite[Statement~4.2]{Veshyn2011}. The proof of this result for matrices which satisfy condition $\Cond$ can be carried out by similar arguments.
\end{proof}

%\subsection{Invertibility for the incompressible vectors} \label{sec: incomp}

For the incompressible vectors, we shall reduce the invertibility problem to a lower bound on the distance between a random vector and a random hyperplane. For this aim we recall Lemma~3.5 from~\cite{RudVesh2008}

\begin{lemma}\label{l:distance}
Let $\A$  be a random matrix from theorem. Let $A_1, ... , A_n$ denote the column vectors of $\A$, and let $H_k$ denote the span of all columns except the $k$-th. Then for every $\delta, \tau \in (0,1)$ and every $\varepsilon > 0$, one has
\begin{equation} \label{eq:r.h.s. distance}
\Pb ( \inf_{x \in Incomp(\delta, \tau)} ||\A x||_2 < \varepsilon n^{-1}) \le \frac{1}{\delta n} \sum_{k=1}^n \Pb (\dist(A_k, H_k) < \tau^{-1} \varepsilon).
\end{equation}
\end{lemma}

%\subsection{Distance via the small ball probability}

Lemma~\ref{l:distance} reduces the invertibility problem to a lower bound on the distance between a random vector and a random hyperplane.

We decompose matrix $\A = \X - z \I$ into the blocks
\begin{equation}\label{eq:matrix}
\begin{pmatrix}
a_{11}  & V^T\\
U   &  \B \\
\end{pmatrix}
\end{equation}
where $\B$ is $ (n-1) \times (n-1)$ matrix, $U, V \in \R^{n-1}$.

Let $h$ be any unit vector orthogonal to $A_2, ... , A_{n}$. It follows that
$$
0 = \begin{pmatrix}
  V^T \\
  \B
\end{pmatrix}^T h = h_1 V + \B^T g,
$$
where $h = (h_1, g)$, and
$$
g = - h_1 \B^{-T} V
$$
From definition of $h$
$$
1 = ||h||_2^2 = |h_1|^2 + ||g||_2^2 = |h_1|^2  + |h_1|^2 ||\B^{-T} V||_2^2
$$
Using this equations we estimate distance
$$
\dist(A_1, H) \geq |(A_1,h)| = \frac{|a_{11} - (\B^{-T} V, U)|}{\sqrt{1+||\B^{-T} V||_2^2}}
$$
It is easy to show that $||\B|| \leq ||\A||$. Let vector $e_1 \in S^{n-2}$ be such that $||\B|| = ||\B e_1 ||_2$. Then we can take vector $e = (0, e_1)^T \in S^{n-1}$ and for this vector
$$
||\A|| \geq ||\A e||_2 = ||(V^T e_1, \B e_1)^T||_2 \geq || \B e_1||_2 = ||\B||.
$$

The bound for right hand sand of~\eqref{eq:r.h.s. distance} will follow from the
\begin{lemma} \label{l:dist vai q.f}
Let matrix $\A$ be from Theorem~\ref{th:smallest singular value}. Then for all $\varepsilon > 0$
\begin{equation} \label{eq:concentration}
   \sup_{v \in \R} \Pb \left (\frac{|(\B^{-T} V, U) - v|}{\sqrt{1+||\B^{-T} V||_2^2}} \le \varepsilon, ||\B|| \leq 3 K \sqrt n \right ) \le C(\rho) \varepsilon^{1/8} + C'(\rho) n^{-1/8},
\end{equation}
where $\B, U, V$ are determined by~\eqref{eq:matrix} and $C(\rho), C_1(\rho)$ are some constants which can depend only on $\rho, K$ and $M_4$.
\end{lemma}

To get this bound we need several statements. We introduce matrix
\begin{equation}
\Q = \begin{pmatrix}
\OO_{n-1}  & \B^{-T}\\
\B^{-1} &  \OO_{n-1}\\
\end{pmatrix}
\quad
W = \begin{pmatrix}
U  \\
V
\end{pmatrix},
\end{equation}
where $\OO_{n-1}$ is $(n - 1) \times (n - 1)$ matrix with zero entries. Scalar product in~\eqref{eq:concentration} can be rewritten using definition of $Q$:
\begin{equation}
\sup_{v \in \R} \Pb \left (\frac{|(\Q W, W) - v|}{\sqrt{1+||\B^{-T}V||_2^2}} \le 2 \varepsilon \right ).
\end{equation}

Introduce vectors
\begin{equation}
W' =\begin{pmatrix}
   U'\\
   V'
\end{pmatrix}
\quad
Z =\begin{pmatrix}
   U\\
   V'
\end{pmatrix},
\end{equation}
where $U', U'$ are independent copies of $U, V$ respectively. We need the following

\begin{statement}\label{l:q.form decoupling}
$$
\sup_{v \in \R} \Pb_W \left (|(\Q W, W) - v| \le 2 \varepsilon \right ) \le \Pb_{W, W'} \left (|(\Q \OP_{J^{c}}(W- W'), \OP_J W) - u| \le 2 \varepsilon \right ),
$$
where $u$ doesn't depend on $\OP_J W = (\OP_J U, \OP_J V)^T$.
\end{statement}
\begin{proof}
Let us fix $v$ and denote
$$
p:=\Pb \left (|(\Q W, W) - v| \le  2 \varepsilon \right ).
$$
We can decompose the set $[n]$ into union $[n] = J \cup J^c$. We can take $U_1 = \OP_J U, U_2 = \OP_{J^c} U, V_1 = \OP_J V$ and $V_2 = \OP_{J^c} V$. By Lemma~\ref{l,a:decoupling}
\begin{align} \label{eq:p}
&p^2 \le \Pb \left (|(\Q W, W) - v| \le  2 \varepsilon , |(\Q   Z,  Z) - v| \le  2 \varepsilon \right )\\
&\le \Pb \left (|(\Q W, W)  - (\Q Z, Z)| \le  4 \varepsilon \right ) \nonumber .
\end{align}
Let us rewrite $\B^{-T}$ in the block form
$$
\B^{-T} = \begin{pmatrix}
\EE  & \F\\
\G &  \HH\\
\end{pmatrix}.
$$
We have
\begin{align*}
& (\Q W, W) = (\EE V_1 , U_1) + (\F V_2, U_1) + (\G V_1, U_2) + (\HH V_2, U_2)\\
& + (\EE^T U_1, V_1) + (\G^T U_2, V_1) + (\F^T U_1, V_2) + (\HH^T U_2, V_2)\\
& (\Q Z, Z) = (\EE V_1, U_1) + (\F V_2', U_1) + (\G V_1, U_2') + (\HH V_2', U_2')  \\
& + (\EE^T U_1, V_1) + (\G^T U_2', V_1) + (\F^T U_1, V_2') + (\HH^T U_2', V_2')
\end{align*}
and
\begin{align}\label{eq:q.f.1}
&(\Q W, W) - (\Q Z,Z) = 2(\F (V_2 - V_2'), U_1) + 2 (\G^T (U_2 - U_2'), V_1)\\
& + 2 (\HH V_2, V_2) - 2(\HH V_2', V_2') \nonumber .
\end{align}
The last two terms in~\eqref{eq:q.f.1} depend only on $U_2, U_2',V_2,V_2'$ and we conclude that
$$
p_1^2 \le \Pb \left ( |(\Q P_{J^c}(W - W'), P_J W ) - u| \le 2 \varepsilon \right),
$$
where $u = u(U_2,V_2,U_2',V_2', \F, \G, \HH)$.
\end{proof}

\begin{statement} \label{l:incomp of inverse}
For all $u \in \R^{n-1}$
$$
\mathbb P\left( \frac{\B^{-T}u}{||\B^{-T}u||_2} \in Comp(\delta, \tau) \text{ and } ||\B|| \le 3 K n^{1/2} \right) \le 2 e^{-c n }.
$$
\end{statement}

\begin{proof}
Let $x = \B^{-T}u$. It is easy to see that
\begin{align*}
\left \{\frac{\B^{-T}u}{||\B^{-T}u||_2} \in Comp(\delta, \tau)\right \} \subseteqq \left \{\exists x: \frac{x}{||x||_2} \in Comp(\delta, \tau) \text{ and } \B^T x = u \right \}
\end{align*}
Replacing matrix $\A$ with $\B^T$ one can easily check that the proof of Lemma~\ref{l:comp. vectors} remains valid for $\B^T$ as well as for $\A$.
\end{proof}

\begin{remark}
  The Statement~\ref{l:incomp of inverse} holds true for $\B^{-T}$ replaced with $\B^{-1}$.
\end{remark}

\begin{statement} \label{l:norm est}
Let $\A$ satisfies condition $\Cond$ and $\B$ be a matrix from decomposition~\eqref{eq:matrix}.
Assume that $||\B|| \le 3 K \sqrt n$. Then with probability at least $1 - e^{-c n}$ matrix $\B$ has the following properties:

\begin{itemize}
\item[a)] $||\B^{-T} V||_2 \geq C$ with probability $1 - e^{-c' n}$ in $W$,\\
\item[b)] $||\B^{-T} V||_2 \leq \varepsilon^{-1/2} ||\B^{-T}||_{HS}$ with probability $1 - \varepsilon$ in $V$,\\
\item[c)] $||\Q W||_2 \geq \varepsilon ||\B^{-T}||_{HS}$ with probability $1 - C'( \varepsilon + n^{-1/2})$ in $W$.
\end{itemize}
\end{statement}

\begin{proof}
Let $\{e_k\}_{k=1}^n$ be a standard basis in $\R^n$. For all $1 \le k \le n$ define vectors by
$$
x_k := \frac{\B^{-1} e_k}{||\B^{-1} e_k||_2}.
$$
By Statement~\ref{l:incomp of inverse} vector $x_k$ is incompressible with probability $1 - e^{-c n}$. We fix matrix $\B$ with such property.\\
a) By norm inequality $||V||_2 \le ||\B|| ||\B^{-T} V||_2$. We know that $ ||\B|| \le 3 K \sqrt n$.
By Lemma~\ref{l,a:rand var levy conc} and Lemma~\ref{l.a:tensorisation} $||V||_2 \geq \sqrt n$.
So we have that $||\B^{-1} V||_2 \geq C$ with probability $1 - e^{-c' n}$. \\
b) By definition
$$
||\B^{-T} V||_2^2 = \sum_{i=1}^n (\B^{-1} e_k, V)^2 = \sum_{i=1}^n ||\B^{-1} e_i||_2^2 (x_k, V)^2.
$$
It is easy to see that $\E(V ,x_k)^2 = 1$. So
$$
\E ||\B^{-T}V||_2^2 = \sum_{i=1}^n ||\B^{-1} e_i||_2^2 = ||\B^{-1}||_{HS}^2.
$$
By Markov inequality
\begin{align*}
&\Pb (||\B^{-T} V||_2 \geq \varepsilon^{-1/2} ||\B^{-1}||_{HS} ) \le \varepsilon.
\end{align*}
c)
By Lemma~\ref{l,a:sum of r.v}, Lemma~\ref{l,a:incomp vec}, Lemma~\ref{l,a:reduction} and Remark~\ref{r:SBP CLT}
\begin{align*}
&\Pb (||\Q W ||_2 \le \varepsilon ||\B^{-1}||_{HS} ) \le \Pb (||\B^{-T} V||_2 \le \varepsilon ||B^{-1}||_{HS} ) \\
&= \Pb (||\B^{-T} V||_2^2 \le \varepsilon ||\B^{-1}||_{HS}^2 ) = \Pb (\sum_{i=1}^n ||\B^{-1} e_i||_2 (x_i, V)^2 \le \varepsilon^2 ||\B^{-1}||_{HS}^2 ) \\
& =\Pb (\sum_{i=1}^n p_i (x_i, V)^2 \le \varepsilon^2  ) \le 2\sum_{i=1}^n p_i \Pb ( (x_i, V) \le \sqrt 2 \varepsilon) \le C' (\varepsilon + n^{-1/2}).
\end{align*}
\end{proof}

\begin{proof}(proof of Lemma~\ref{l:dist vai q.f})
Let $\xi_1, ... ,\xi_n$ be i.i.d. Bernoulli random variables  with $\E \xi_i = c_{0}/2$. We define $J: =\{i: \xi_i = 0 \}$ and $\E_0 : = \{|J^c| \le c_0 n \}$.
From large deviation inequality we may conclude that $\Pb(E_0) \geq 1 - 2 \exp(-c_0^2 n/2)$. Introduce event
$$
E_1:= \{\varepsilon_0^{1/2} \sqrt {1 + ||\B^{-T} V||_2^2} \le ||\B^{-1}||_{HS} \le \varepsilon_0^{-1}||\Q \OP_{J^c} (W - W')||_2 \},
$$
where $\varepsilon_0$ will be choosen later.

From Statement~\ref{l:norm est} we can conclude that
$$
\Pb_{\B,W,W', J} (E_1 \cup ||B|| \geq 3 K \sqrt n) \geq 1 - C' ( \varepsilon_0 + n^{-1/2}) - 2 e^{-c' n}.
$$
Consider the random vector
$$
w_0 = \frac{1}{||\Q \OP_{J^c} (W - W')||_2}\begin{pmatrix}
  \B^{-T} \OP_{J^c} (V - V') \\
  \B^{-1} \OP_{J^c} (U - U')
\end{pmatrix} =
\begin{pmatrix}
   a\\
   b
\end{pmatrix}.
$$
By Statement~\ref{l:incomp of inverse} it follows that the event $E_2:=\{a \in incomp(\delta,\tau) \}$ holds with probability
$$
\Pb_{\B}(E_2 \cup ||\B|| \geq 3 K \sqrt n | W, W', J) \geq 1 -  2 \exp(-c'' n).
$$
Combining these probabilities we have
\begin{align*}
& \Pb_{\B,W,W', J} (E_0,~E_1,~E_2 \cup ||\B|| \geq 3 K \sqrt n ) \\
& \geq 1 -  2 e^{-c_{0}^2 n/2 } - C'( \varepsilon_0 + n^{-1/2}) - 2 e^{-c' n} - 2 e^{-c'' n}: = 1 - p_0.
\end{align*}
We may fix $J$ that satisfies $|J^c| \le c_0$  and
$$
\Pb_{\B,W,W'} (E_1,~E_2 \cup ||\B|| \geq 3 K \sqrt n ) \geq 1 - p_0.
$$
By Fubini's theorem $\B$ has the following property with probability at least $1 - \sqrt{p_0}$
$$
\Pb_{W,W'} (E_1,~E_2 \cup ||\B|| \geq 3 K \sqrt n | \B ) \geq 1 - \sqrt{p_0}.
$$
The event $\{||\B|| \geq 3 K \sqrt n\}$ depends only on $\B$. We may conclude that random matrix $\B$ has the following property with probability at least $1 - \sqrt{p_0}$: either $||\B|| \geq 3 K \sqrt n$, or
\begin{align} \label{eq:B prop}
||\B|| \le 3 K \sqrt n \text{ and } \Pb_{W,W'} (E_1,~E_2| \B) \geq 1 - \sqrt p_0
\end{align}
The event we are interested in is
$$
\Omega_0 := \left ( \frac{|(\Q W, W) - u|}{\sqrt{1 + ||\B^{-T} V||_2^2}} \le 2 \varepsilon \right ).
$$
We need to estimate probability
$$
\Pb_{\B, W} (\Omega_0 \cap ||\B|| \le 3 K \sqrt n) \le \Pb_{\B,W} (\Omega_0 \cap \text{ \eqref{eq:B prop} holds} ) + \Pb_{\B,W} (\Omega_0 \cap \text{ \eqref{eq:B prop} fails} ).
$$
The last term is bounded by $\sqrt{p_0}$.
$$
\Pb_{\B,W} (\Omega_0 \cap ||\B|| \le 3 K \sqrt n) \le \sup_{\B \text{ satisfies~\eqref{eq:B prop}}} \Pb_{W} (\Omega_0| \B ) + \sqrt{p_0}.
$$
We can conclude that
\begin{align*}
\Pb_{\B,W} (\Omega_0 \cap ||\B|| \le 3 K \sqrt n) \le \sup_{\B \text{ satisfies~\eqref{eq:B prop}}} \Pb_{W,W'} (\Omega_0,~E_1| \B ) + 2\sqrt{p_0}.
\end{align*}
Let us fix $\B$ that satisfies~\eqref{eq:B prop} and denote $p_1: = \Pb_{W,W'} (\Omega_0,~E_1| \B )$. By Statement~\ref{l:q.form decoupling} and the first inequality in $E_1$ we have
\begin{align*}
p_1^2 \le \Pb_{W, W'} \left ( \underbrace{|(\Q \OP_{J^c} (W - W'), \OP_J W) - v| \le \frac{\varepsilon}{\sqrt{\varepsilon_0} } ||\B^{-1}||_{HS}}_{\Omega_1} \right )
\end{align*}
and
\begin{align*}
\Pb_{W, W'} ( \Omega_1 ) \le \Pb_{W, W'} ( \Omega_1,~E_1,~E_2) +\sqrt{p_0}.
\end{align*}
Further
$$
p_1^2 \le \Pb_{W,W'} ( |(w_0, \OP_J W) - v| \le 2\varepsilon_0^{-3/2} \varepsilon,~E_2 ) + \sqrt{p_0}.
$$
By definition random vector $w_0$ is determined by the random vector $\OP_{J^c}(W - W')$, which is independent of the random vector $\OP_J W$. We fix $\OP_{J^c}(W - W')$ and have
$$
p_1^2 \le \sup_{\substack {w_0 = (a,b)^{T}:\\ a \in Incomp(\delta,\tau)\\ w \in \R}} \Pb_{\OP_J W} \left ( |(w_0, \OP_J W) - w| \le \varepsilon_0^{-3/2} \varepsilon \right ) + \sqrt{p_0}.
$$
Let us fix a vector $w_0$ and a number $w$. We can rewrite
\begin{equation}\label{eq:sum}
(w_0, P_J W)  = \sum_{i \in J} ( a_i x_i + b_i y_i),
\end{equation}
where $||a||_2^2 + ||b||_2^2 = 1$. From Lemma~\ref{l,a:incomp vec} and Remark~\ref{r:incomp vec} we know that at least $[2 c_{0} n]$ coordinates of vector $a \in Incomp(\delta,\tau)$ satisfy
$$
\frac{\tau}{\sqrt{2 n}} \le |a_k| \le \frac{1}{\sqrt{\delta n}}.
$$
We denote the set of coordinates of $a$ with this property by $\spr(a)$. By construction of $J$ we can conclude that $|\spr(a)| = [c_{0} n]$.
By Lemma~\ref{l,a:reduction} we can reduce our sum~\eqref{eq:sum} to the set $\spr(a)$.
Now we will find the properties of $|b_i|$. We can decompose the set $\spr(a)$ into two sets $\spr(a) = I_1 \cup I_2$: \\
a) $I_1 = \{ i \in \spr(a): |b_i| \sqrt n \rightarrow \infty $ as $n \rightarrow \infty \}$; \\
c) $I_2 = \{i \in \spr(a):  |b_i| = O(n^{-1/2})\}$;\\
From $||b||_2^2 < 1$ it follows that $|I_1| = o(n)$. For $i \in I_2$ we have $|a_i^{-1} b_i|  =  O(1)$. By Lemma~\ref{l,a:reduction} we have
$$
\Pb (|\sum_{i \in \spr(a)} (a_i x_i + b_i y_i) - w| < 2 \varepsilon_0^{-3/2} \varepsilon) \le \Pb (|\sum_{i \in I_2} (a_i x_i + b_i y_i) - w'| < 2 \varepsilon_0^{-3/2} \varepsilon).
$$
We can apply Statement~\ref{st:SBP CLT}
$$
\Pb (|\sum_{i \in I_2} (a_i x_i + b_i y_i) - w'| < 2 \varepsilon_0^{-3/2} \varepsilon) \le \frac{C_1 \varepsilon_0^{-3/2} \varepsilon}{(1-\rho^2)^{1/2}}  + C_2 (1-\rho^2)^{-3/2} n^{-1/2}.
$$
It follows that
\begin{align*}
&\Pb_{\B, W} (\Omega_0 \cap ||\B|| \le 3 K \sqrt n) \le \\
&\le\left (\frac{C_1 \varepsilon_0^{-3/2} \varepsilon}{(1-\rho^2)^{1/2}}  + C_2 (1-\rho^2)^{-3/2} n^{-1/2} \right)^{1/2} + p_0^{1/4} + 2 \sqrt{p_0}.
\end{align*}
We take $\varepsilon_0 = \varepsilon^{1/2}$ and conclude that
$$
\Pb_{\B, W} (\Omega_0 \cap ||\B|| \le 3 K \sqrt n) \le C(\rho) \varepsilon^{1/8} + C'(\rho) n^{-1/8},
$$
where $C(\rho), C'(\rho)$ are some constants which depend on $\rho, K$ and $M_4$.
\end{proof}

\begin{proof}(proof of Theorem~\ref{th:smallest singular value})
The result of the theorem follows from Lemmas~\ref{l:comp. vectors}, \ref{l:distance} and~\ref{l:dist vai q.f}.
\end{proof}

\begin{remark} \label{r: smallest sing value}
It not very difficult to show that we can change matrix $z \I$ in Theorem~\ref{th:smallest singular value}
by arbitrary non-random matrix $\M$ with $||\M|| \le K \sqrt n$. We can also assume that $\E X_{ij}^2 \geq 1$.
Results of section~\ref{sec: decomposition} are based on Lemmas~\ref{l.a:sum. c.f.} and~\ref{l.a:tensorisation} which doesn't depend on shifts. It is easy to see
that Statement~\ref{l:norm est} still holds true if we assume that $\varepsilon < n^{-Q}$ for some $Q > 0$.
Then we can reformulate Theorem~\ref{th:smallest singular value} in the following way: there exist some constants $A,B > 0$ such that
$$
\Pb (s_n(\X + \M) \le \varepsilon n^{-A}, ||\X + \M|| \leq K \sqrt n) \le C(\rho) n^{-B}.
$$
\end{remark}

\section{Uniform integrability of logarithm} \label{sec:uniform}

In this section we prove the next result

\begin{theorem}\label{th: log uniform integr}
Under the condition $\Cond$ $\log( \cdot )$ is uniformly integrable in probability with respect to $\{\nu_n\}_{n \geq 1}.$
\end{theorem}

Before we need several lemmas about the behavior of the singular values

\begin{lemma} \label{l:largest singular value}
If condition $\Cond$ holds then there exists constant $K: = K(\rho)$ such that $\Pb(s_1(\X) \geq K \sqrt n) = o(1)$.
\end{lemma}
\begin{proof}
We can decompose matrix $\X$ into symmetric and skew-symmetric matrices:
$$
\X = \frac{\X + \X^T}{2} + \frac{\X - \X^T}{2}  = \X_1 + \X_2.
$$
In~\cite[Theorem~2.3.23]{Tao2012} it is proved that for some $K_1 > \sqrt{2(1 + \rho)}$
\begin{equation}\label{eq: s x1}
\Pb(s_1( \X_1) \geq K_1 \sqrt n) = o(1).
\end{equation}
and for some $K_2 > \sqrt{2(1 - \rho)}$
\begin{equation}\label{eq: s x2}
\Pb(s_1(\X_2) \geq K_2 \sqrt n) = o(1)
\end{equation}
Set $K = 2\max(K_1, K_2)$. From~\eqref{eq: s x1},~\eqref{eq: s x2} and inequality
$$
s_1(\X) \le s_1(\X_1) + s_1(\X_2)
$$
it follows that
\begin{align*}
&\Pb(s_1(\X) \geq K \sqrt n) \le \Pb \left (\left \{s_1(\X_1) \geq \frac{K \sqrt n }{2}\right \} \cup \left \{ s_1(\X_2) \geq \frac{K \sqrt n }{2} \right \}\right ) \\
&\le \Pb \left  (s_1( \X_1) \geq \frac{K \sqrt n }{2}\right ) + \Pb \left  (s_1(\X_2) \geq \frac{K \sqrt n }{2} \right ) = o(1).
\end{align*}
\end{proof}

\begin{remark} \label{r:largest singular value}
Suppose that elements of $\X_n$ depend on $n$, but satisfy conditions $\Cond$ and $|x_{ij}| \le \delta_n \sqrt n, \E x_{ij}^2 \le 1$
and $\E|x_{ij}|^l \le b(\delta_n \sqrt n)^{l-1}$ for some $b>0, l \geq 3$ and $\delta_n \rightarrow 0$
with the convergence rate slower that any preassigned one as $n \rightarrow \infty$.
Then for some $K > 0$ it can be shown that
$$
\Pb(s_1(\X) \geq K \sqrt n) = o(n^{-l}).
$$
For the proof see~\cite[Theorem~5.1]{BaiSilv2010}.
\end{remark}

\begin{lemma}\label{l:large singular values}
If condition $\Cond$ holds then there exist $c > 0$ and $0 < \gamma < 1$ such that a.s. for $n \gg 1$ and $n^{1-\gamma} \le i \le n-1$
$$
s_{n-i}(n^{-1/2}\X - z\I) \geq c \frac{i}{n}.
$$
\end{lemma}
\begin{proof}
Set $s_i: = s_i(n^{-1/2}\X - z\I)$. Up to increasing $\gamma$, it is sufficient to prove
the statement for all $2 (n-1)^{1-\gamma} \le i \le n-1$ for some $\gamma \in (0,1)$ to be chosen later.
We fix some $2 (n-1)^{1-\gamma} \le i \le n-1$ and consider the matrix $\A'$ formed by
the first $m:=n- \lceil i/2 \rceil$ rows of $\sqrt n \A$. Let $s_1' \geq ... \geq s_m'$ be the singular values of $\A'$. We get
$$
n^{-1/2} s_{n-i}' \le s_{n-i}.
$$
By $R_i$ we denote the row of $\A'$ and $H_i = \Sp(R_j, j = 1, ..., m, j \neq i)$. By Lemma~\ref{l,a:distance} we obtain
$$
s_1'^{-2} + ... +s_{n- \lceil i/2 \rceil}'^{-2} = \dist_1^{-2} + ... + \dist_{n- \lceil i/2 \rceil}^{-2}.
$$
We have
\begin{align}\label{eq: chain}
\frac{i}{2n} s_{n-i}^{-2} \le \frac{i}{2} s_{n-i}'^{-2} \le \sum_{j= n - i}^{ n - \lceil i/2 \rceil } s_j'^{-2} \le \sum_{j=1}^{n- \lceil i/2  \rceil} \dist_j^{-2},
\end{align}
where $\dist_j := \dist(R_j,H_j)$. To estimate $\dist(R_j, H_j)$ we would like to apply Lemma~\ref{l,a:distance est},
but we can't do it directly, because $R_j$ and $H_j$ are not independent. Let's consider the case $j = 1$ only.
To estimate distance $\dist_1$ we decompose matrix $\A'$ into the blocks
$$
\A' = \begin{pmatrix}
a_{1,1}  & Y\\
X &  \B
\end{pmatrix},
$$
where $X \in \R^{m-1}, Y^{T} \in \R^{n-1}$ and $\B$ is an ${m-1} \times n-1$ matrix formed by
rows $B_1, ..., B_{m-1}$. We denote by $H_1' = \Sp(B_1, ..., B_{m-1})$. From definition of distance
$$
\dist(R_1, H_1) = \inf_{v \in H_1} ||R_1 - v||_2 \geq \inf_{u \in H'} ||Y - u||_2 = \dist(Y, H_1')
$$
and
$$
\dim(H_1') \le \dim(H_1) \le n - 1 - i/2 \le n - 1 - (n - 1)^{1 - \gamma}.
$$
Now vector $Y$ and hyperplane $H_1'$ are independent. Fixing realization of $H_1'$, by Lemma~\ref{l,a:distance est}, with $n, R, H$ replaced with $ n - 1, Y, H_1'$
respectively, we can obtain that
$$
\mathbb P( \dist(Y,H_1') \le \frac{1}{2} \sqrt{n-1 - \dim(H_1')}) \le \exp(-(n-1)^\delta).
$$
Using this inequality it is easy to show that
$$
\Pb \left ( \bigcup_{n \gg 1} \bigcup_{i= \lceil 2 (n-1)^{1-\gamma} \rceil}^{n-1} \bigcup_{j=1}^{n-\lceil i/2 \rceil} \left \{\dist(R_j,H_j) \le \frac{1}{2} \sqrt{\frac{i}{2}} \right \} \right ) < \infty.
$$
Now by Borel-Cantelli lemma and~\eqref{eq: chain}  we can conclude the statement of the lemma.
\end{proof}
\begin{remark}\label{r:large singular values}
Lemma~\ref{l:large singular values}  holds true if we assume that $\E X_{ij} \neq 0$ and $\E X_{i j}^2 = 1 + o(1)$.
\end{remark}

\begin{proof} (Proof of Theorem~\ref{th: log uniform integr})
To prove Theorem~\ref{th: log uniform integr} we need to show that there exist $p, q > 0$ such that
\begin{equation} \label{eq:log+ u.i.}
\lim_{t \rightarrow \infty} \varlimsup_{n \rightarrow \infty} \mathbb P \left (\int_0^\infty x^{p} \nu_n(dx) > t \right ) = 0
\end{equation}
and
\begin{equation} \label{eq:log- u.i.}
\lim_{t \rightarrow \infty} \varlimsup_{n \rightarrow \infty} \mathbb P \left (\int_0^\infty x^{-q} \nu_n(dx) > t \right ) = 0.
\end{equation}

By Lemma~\ref{l:largest singular value} there exists set $\Omega_0:=\Omega_{0, n} = \{\omega \in \Omega: s_1(\X) \le K n^{1/2} \}$ such that
\begin{equation} \label{eq: largest s.v}
\Pb(\Omega_{0}) = 1 - o(1).
\end{equation}
We conclude~\eqref{eq:log+ u.i.} from~\eqref{eq: largest s.v} for $p = 2$.

We denote $\Omega_1:= \Omega_{1,n} = \{\omega \in \Omega: s_{n-i} > c \frac{i}{n}, n^{1-\gamma} \le i \le n-1\}$.
Let us consider the set $\Omega_2:=\Omega_{2,n} = \Omega_{1} \cap \{\omega: s_n \geq n^{-B-1/2}\}$, where $B > 0$.
We decompose probability from~\eqref{eq:log- u.i.} into two terms
$$
\Pb \left (\int_0^\infty x^{-q} \nu_n(dx) > t \right ) = \mathbb I_1 + \mathbb I_2,
$$
where
\begin{align*}
&\mathbb I_1: = \Pb \left (\int_0^\infty x^{-q} \nu_n(dx) > t, \Omega_{2} \right),\\
&\mathbb I_2: = \Pb \left (\int_0^\infty x^{-q} \nu_n(dx) > t,  \Omega_{2}^c \right).
\end{align*}
We can estimate $\mathbb I_2$ by
$$
\mathbb I_2 \le \Pb(s_n(\X - \sqrt n z \I) \le n^{-A}, \Omega_{0}) + \Pb(\Omega_{0}^c) + \Pb(\Omega_{1}^c).
$$
From Theorem~\ref{th:smallest singular value} it follows that
\begin{equation}\label{eq: least sv}
\Pb(s_n(\X - \sqrt n z \I) \le n^{-B}, \Omega_{0}) \le C(\rho)n^{-1/8}.
\end{equation}
By Lemma~\ref{l:large singular values}
\begin{equation} \label{eq: large sv}
\varlimsup_{n \rightarrow \infty} \Pb(\Omega_{1}^c) = 0.
\end{equation}
From~\eqref{eq: largest s.v},~\eqref{eq: least sv} and~\eqref{eq: large sv} we conclude
$$
\varlimsup_{n \rightarrow \infty} \mathbb I_2 = 0.
$$
To prove~\eqref{eq:log- u.i.} it remains to bound  $\mathbb I_1$. From Markov inequality
$$
\mathbb I_1 \le \frac{1}{t} \E \left[\int_0^\infty x^{-q} \nu_n(dx) \one(\Omega_{2}) \right ].
$$
By definition of $\Omega_{2}$
\begin{align*}
\E \left[\int x^{-q} \nu_n(dx) \one(\Omega_2) \right ] \le \frac{1}{n} \sum_{i=1}^{n-\lceil n^{1-\gamma} \rceil} s_{i}^{-q} + \frac{1}{n}\sum_{i=n - \lceil n^{1-\gamma} \rceil + 1}^{n} s_{i}^{-q} \\
\le 2 n^{q (B  + 1/2) - \gamma} + c^{-q} \frac{1}{n} \sum_{i=1}^n \left ( \frac{n}{i}\right )^q \le 2 n^{q (B  + 1/2) - \gamma} + c^{-q}\int_{0}^1 s^{-q} ds.
\end{align*}
If $0 < q < \min(1, \gamma/(B + 1/2))$ then the last integral is finite.
\end{proof}

\section{Convergence of singular values} \label{sec:singular values}

Let function $\mathcal{F}_n(x,z)$ be an empirical distribution function of singular values $s_1 \geq ... \geq s_n$
of matrix $n^{-1/2}\X - z\I$ which corresponds to measure $\nu_n(z, \cdot)$.

Let us recall definition of Stieltjes transform
\begin{definition}
The Stieltjes transform of measure $m(\cdot)$ on $\R$ is
$$
S(\alpha) = \int_{\R} \frac{m(dx)}{x - \alpha}, \quad \alpha \in \C^{+}.
$$
\end{definition}

In this section we prove the following theorem
\begin{theorem}\label{th:Stieltjes transform}
Assume that condition $\Cond$ holds true. There exists non-random distribution function $\mathcal{F}(x,z)$ such that
for all continues and bounded functions $f(x)$, a.a. $z \in \C$ and all $\varepsilon > 0$
\begin{equation*}
\Pb \left ( \left |\int_\R f(x) d\mathcal{F}_n(x,z) - \int_\R f(x) d\mathcal{F}(x,z) \right | > \varepsilon \right ) \rightarrow 0 \text{ as } n \rightarrow \infty,
\end{equation*}
\end{theorem}
\begin{proof}

First we show that family $\{\mathcal{F}(z, x)\}_{n \geq 1}$ is tight. From strong law of large numbers it follows that
$$
\int_0^\infty x^2 d\mathcal{F}(x,z) \le \frac{1}{n^2} \sum_{i,j=1}^n X_{ij}^2 \rightarrow 1 \text{ as } n \rightarrow \infty.
$$
Using this and the fact that $s_i(n^{-1/2}\X - z \I) \le s_i(n^{-1/2}\X) + |z|$ we conclude tightness of  $\{\mathcal{F}_n(z, x)\}_{n \geq 1}$.
If we show that $\mathcal{F}_n$ weakly converges in probability to some function $\mathcal{F}$, then $\mathcal{F}$ will be distribution function.

Introduce the following $2 n \times 2n$ matrices
\begin{equation} \label{eq:matrices}
\V = \begin{pmatrix}
\OO_n  &  n^{-1/2}\X\\
 n^{-1/2} \X^{T} & \OO_n\\
\end{pmatrix}, \quad
\J(z) = \begin{pmatrix}
\OO_n  &  z \I\\
\overline z \I &  \OO_n\\
\end{pmatrix}
\end{equation}
where $\OO_n$ denotes $n \times n$ matrix with zero entries. Consider matrix
$$
\V(z) := \V -\J(z).
$$
It is known that eigenvalues of $\V(z)$ are singular values of $n^{-1/2} \X - z \I$ with signs $\pm$.

It is easy to see that empirical distribution function $F_n(x,z)$ of eigenvalues of matrix $\V(z)$ can be written in the following way
$$
F_n(x,z) = \frac{1}{2n} \sum_{i=1}^n \one\{s_i \le x\} + \frac{1}{2n} \sum_{i=1}^n \one\{-s_i \le x\}.
$$
There is one to one correspondence between $\mathcal{F}_n(x,z)$ and $F_n(x,z)$
$$
F_n(x,z) = \frac{1 + \sgn(x) \mathcal{F}_n(|x|,z)}{2}
$$
So it is enough to show that there exists non-random distribution function $F(x,z)$ such that for all continues and bounded functions $f(x)$, and a.a. $z \in \C$
\begin{equation} \label{eq: sv convergence}
\Pb \left ( \left |\int_\R f(x) d{F}_n(x,z) - \int_\R f(x) d{F}(x,z) \right | > \varepsilon \right ) \rightarrow 0 \text{ as } n \rightarrow \infty.
\end{equation}

We denote Stieltjes transforms of $F_n$ and $F$ by $S_n(x, z)$ and $S(x, z)$ respectively.
Due to the relations between distribution functions and Stieltjes transforms,~\eqref{eq: sv convergence}
will follow from
\begin{equation} \label{eq: sv convergence2}
\Pb(|S_n(\alpha, z) - S(\alpha, z)| > \varepsilon) \rightarrow 0 \text{ as } n \rightarrow \infty,
\end{equation}
for a.a. $z \in \C$ and all $\alpha \in \C^+$.

Set
\begin{equation} \label{eq:resolvent}
\RR(\alpha, z):= (\V(z) - \alpha \I_{2n})^{-1}.
\end{equation}
By definition $S_n(\alpha,z ) = \frac{1}{2n} \Tr \RR(\alpha, z)$. We introduce the following function
$$
s_n(\alpha, z) := \E S_n(\alpha,z ) = \frac{1}{2 n} \sum_{i=1}^{2 n} \E [\RR(\alpha, z)]_{i i},
$$
One can show that
$$
s_n(\alpha, z) = \frac{1}{n} \sum_{i=1}^{n} \E [\RR(\alpha, z)]_{i i} = \frac{1}{n} \sum_{i=n+1}^{2 n} \E [\RR(\alpha, z)]_{i i}
$$
By Chebyshev inequality and Lemma~\ref{stilt.transform} it is staighforward to check that
\begin{equation} \label{eq: sv convergenes3}
|s_n(\alpha, z) -  s(\alpha, z)| \rightarrow 0 \text{ as } n \rightarrow \infty.
\end{equation}
implies~\eqref{eq: sv convergence2}.

%Introduce the following $2 n \times 2n$ matrices
%\begin{equation} \label{eq:matrices}
%\V = \begin{pmatrix}
%n^{-1/2}\X  &  \OO_n\\
%\OO_n &  n^{-1/2} \X^{T}\\
%\end{pmatrix}, \quad
%\J(z) = \begin{pmatrix}
%\OO_n  &  z \I\\
%\overline z \I &  \OO_n\\
%\end{pmatrix}, \quad
%\J := \J(1),
%\end{equation}
%where $\OO_n$ denotes $n \times n$ matrix with zero entries. Consider matrix
%$$
%\V(z) := \V\J -\J(z).
%$$
%It is known that eigenvalues of $\V(z)$ are singular values of $n^{-1/2} \X - z \I$ with signs $\pm$.
%
%Set for $\alpha = u + i v, v > 0$
%\begin{equation} \label{eq:resolvent}
%\RR(\alpha, z):= (\V(z) - \alpha \I_{2n})^{-1}.
%\end{equation}
%We introduce the following functions
%$$
%s_n(\alpha, z) = \frac{1}{n} \sum_{i=1}^n \E [\RR(\alpha, z)]_{i i} = \frac{1}{n} \sum_{i=1}^n \E [\RR(\alpha, z)]_{i+n, i+n} = \frac{1}{2 n} \sum_{i=1}^{2 n} \E [\RR(\alpha, z)]_{i i},
%$$
%\begin{equation*}
%t_n(\alpha, z) = \frac{1}{n} \sum_{i=1}^n \E [\RR(\alpha, z)]_{i+n, i}, \quad u_n(\alpha, z) = \frac{1}{n} \sum_{i=1}^n \E [\RR(\alpha, z)]_{i, i+n}.
%\end{equation*}
%It is easy to see that $s_n$ is a Stieltjes transform of
%$$
%\E \hat \nu_{n^{-1/2} \X - z \I}(\cdot) = \frac{\E\nu_n(\cdot) + \E\nu_n(- \cdot)}{2},
%$$
%which is the symmetrized version of a measure $\E \nu_n$. There is a relation between Stieltjes transform $S_n(\alpha, z)$ and $s_n(\alpha,z)$:
%$$
%s_n(\alpha, z) = \alpha S_n(\alpha^2,z).
%$$
%To prove Theorem~\ref{th:Stieltjes transform} we need to show convergence of $s_n(\alpha,z)$ to Stieltjes transform of some symmetric measure $\nu_0(z, \cdot)$.

By resolvent equality we may write
$$
1 + \alpha s_n(\alpha, z) = \frac{1}{2n} \E \Tr (\V \RR(\alpha, z)) - z t_n(\alpha, z) - \overline z u_n(\alpha, z).
$$
Introduce the  notation
$$
\mathbb A := \frac{1}{2 n} \E \Tr (\V \RR)
$$
and represent $\mathbb A$ as follows
$$
\mathbb A = \frac{1}{2} \mathbb A_1 + \frac{1}{2} \mathbb A_2,
$$
where
\begin{equation*}
\mathbb A_1 = \frac{1}{n} \sum_{i=1}^{n} \E [\V  \RR]_{i i}, \quad \mathbb A_2 = \frac{1}{n} \sum_{i=1}^{n} \E [\V \RR]_{i+n, i+n}.
\end{equation*}

First we consider $\mathbb A_1$. By definition of the matrix $\V$, we have
$$
\mathbb A_1 = \frac{1}{n^{3/2}} \sum_{j, k = 1}^n \E X_{j k} R_{ k + n, j}.
$$
Note that
\begin{align*}
\frac{\partial \RR}{\partial X_{j k}} = - \frac{1}{\sqrt n} \RR [e_j e_{k+n}^{T}] \RR.
\end{align*}
Applying Lemma~\ref{l,a:stein} we obtain
$$
\mathbb A_1 = \mathbb B_1 + \mathbb B_2 + \mathbb B_3 + \mathbb B_4  + r_n (\alpha, z).
$$
where
\begin{align*}
& \mathbb B_1 = -\frac{1}{n^2} \sum_{j, k = 1}^ n  \E[\RR [e_j e_{k+n}^{T}] \RR]_{k+n, j} =
-\frac{1}{n^2} \sum_{j, k = 1}^ n  \E (R_{k+n, j})^2 \\
& \mathbb B_2 = -\frac{1}{n^2} \sum_{j, k = 1}^ n  \E[\RR [e_{k+n} e_{j}^{T}] \RR]_{k+n, j} =
-\frac{1}{n^2} \sum_{j, k = 1}^ n  \E R_{j j} R_{k+n, k+n}\\
& \mathbb B_3 = -\frac{\rho}{n^2} \sum_{j, k = 1}^ n  \E[\RR [e_k e_{j+n}^{T}] \RR]_{k+n, j}=
-\frac{\rho}{n^2} \sum_{j, k = 1}^ n  \E R_{k+n, k} R_{j+n,j}\\
& \mathbb B_4 = -\frac{\rho}{n^2} \sum_{j, k = 1}^ n \E[\RR [e_{j+n} e_{k}^{T}] \RR]_{k+n, j} =
-\frac{\rho}{n^2} \sum_{j, k = 1}^ n  \E R_{k j} R_{k+n, j+n}.
\end{align*}

Without loss of generality we can assume further that $\E X_{1 1}^2 = 1$ because the impact of diagonal is of order $O(n^{-1})$.

From $||\RR||_{HS} \le \sqrt{n} ||\RR|| \le \sqrt{n} v^{-1}$ it follows
$$
|\mathbb B_1| \le \frac{1}{n^2} \sum_{j, k = 1}^ n \E X_{j k}^2 \E (R_{k+n, j})^2 \le \frac{1}{n v^2}.
$$
Similarly
$$
|\mathbb B_4| \le \frac{1}{v^2 n}.
$$
By Lemma~\ref{stilt.transform} $\mathbb B_2 = - s_n^2(\alpha, z) + \varepsilon(\alpha, z)$. By Lemma~\ref{l,a:second diag} $\mathbb B_3 = -\rho t_n^2(\alpha, z) + \varepsilon(\alpha, z)$. We obtain that
$$
\mathbb A_1 = -s_n^2(\alpha,z) - \rho t_n^2(\alpha, z) + \delta_n(\alpha, z).
$$

Now we consider the term $\mathbb A_2$. By definition of the matrix $\V$, we have
$$
\mathbb A_2 = \frac{1}{n^{3/2}} \sum_{j, k = 1}^n \E X_{j k} R_{ j,k+n}.
$$
By Lemma~\ref{l,a:stein} we obtain that
\begin{equation}
\mathbb A_2 = \mathbb C_1 + \mathbb C_2 + \mathbb C_3 + \mathbb C_4  + r_n (\alpha, z).
\end{equation}
where
\begin{align*}
& \mathbb C_1 = -\frac{1}{n^2} \sum_{j, k = 1}^ n  \E[\RR [e_j e_{k+n}^{T}] \RR]_{j, k+n} =
-\frac{1}{n^2} \sum_{j, k = 1}^ n \E R_{j j} R_{k+n, k+n}\\
&\mathbb C_2 = -\frac{1}{n^2} \sum_{j, k = 1}^ n  \E[\RR [e_{k+n} e_{j}^{T}] \RR]_{j, k+n} =
-\frac{1}{n^2} \sum_{j, k = 1}^ n  \E (R_{j, k+n})^2\\
&\mathbb C_3 = -\frac{\rho}{n^2} \sum_{j, k = 1}^ n  \E[\RR [e_k e_{j+n}^{T}] \RR]_{j, k+n}=
-\frac{\rho}{n^2} \sum_{j, k = 1}^ n  \E R_{j k} R_{j+n,k+n}\\
&\mathbb C_4 = -\frac{\rho}{n^2} \sum_{j, k = 1}^ n  \E[\RR [e_{j+n} e_{k}^{T}] \RR]_{j, k+n} =
-\frac{\rho}{n^2} \sum_{j, k = 1}^ n  \E R_{j, j+n} R_{k,k+n }.
\end{align*}
It is easy to show that
$$
|\mathbb C_2| \le \frac{1}{v^2 n}, \quad |\mathbb C_3| \le \frac{1}{v^2 n}.
$$
By Lemma~\ref{stilt.transform} $\mathbb C_1 = -s_n^2(\alpha, z) + \varepsilon_n(\alpha, z)$. By Lemma~\ref{l,a:second diag} $\mathbb C_4 = - \rho u_n^2(\alpha, z) + \varepsilon_n(\alpha, z)$. We obtain that
\begin{equation*}
\mathbb A_2 = -s_n^2(\alpha,z) - \rho  u_n^2(\alpha, z) + \delta_n(\alpha, z).
\end{equation*}
So we have that
\begin{equation*}
\mathbb A = -s_n^2(\alpha, z) - \frac{\rho}{2} t_n^2(\alpha,z) - \frac{\rho}{2} u_n^2(\alpha,z)  + \varepsilon_n(\alpha,z).
\end{equation*}

No we will investigate the term $z t_n(\alpha, z)$ which we may represent as follows
$$
\alpha t_n(\alpha, z) = \frac{1}{n} \sum_{j=1}^n \E [\V(z)\RR]_{j+n,j} = \frac{1}{n} \sum_{j=1}^n \mathbb E [\V \RR]_{j+n,j} - \overline z s_n(\alpha, z).
$$
By definition of the matrix $\V$, we have
\begin{align*}
&\alpha t_n(\alpha, z) = \frac{1}{n^{3/2}} \sum_{j, k=1}^n \E  X_{j k} R_{j,k} -  \overline z s_n(\alpha, z) = \\
&\mathbb D_1 + \mathbb D_2 + \mathbb D_3 + \mathbb D_4 - \overline z s_n(\alpha, z) + r_n(\alpha, z),
\end{align*}
where
\begin{align*}
& \mathbb D_1 = -\frac{1}{n^2} \sum_{j, k = 1}^n  \E[\RR [e_j e_{k+n}^{T}] \RR]_{j, k} =
-\frac{1}{n^2} \sum_{j, k = 1}^n  \E R_{j, j} R_{k+n, k} \\
& \mathbb D_2 = -\frac{1}{n^2} \sum_{j, k = 1}^ n  \E[\RR [e_{k+n} e_{j}^{T}] \RR]_{j, k} =
-\frac{1}{n^2} \sum_{j, k = 1}^n  \E R_{j, k+n} R_{j,k}\\
& \mathbb D_3 = -\frac{\rho}{n^2} \sum_{j, k = 1}^ n  \E[\RR [e_k e_{j+n}^{T}] \RR]_{j, k}=
-\frac{\rho}{n^2} \sum_{j, k = 1}^n \E R_{j, k} R_{j+n,k}\\
& \mathbb D_4 = -\frac{\rho}{n^2} \sum_{j, k = 1}^n \E[\RR [e_{j+n} e_{k}^{T}] \RR]_{j, k} =
 -\frac{\rho}{n^2} \sum_{j, k = 1}^n \E R_{j, j+n} R_{k,k}.
\end{align*}
By similar arguments as before we can prove that
$$
|\mathbb D_2| \le \frac{1}{v^2 n}, \quad |\mathbb D_3| \le \frac{1}{v^2 n}
$$
and $\mathbb D_1 = -s_n(\alpha, z) t_n(\alpha, z)+ \varepsilon_n(\alpha, z)$,  $\mathbb D_4 =  -\rho s_n(\alpha, z) u_n(\alpha, z) + \varepsilon_n(\alpha, z)$. We obtain that
$$
\alpha t_n(\alpha, z) = - s_n(\alpha,z) t_n(\alpha, z) - \rho  s_n(\alpha, z) u_n(\alpha, z) - \overline z s_n(\alpha, z) + \delta_n(\alpha, z).
$$

Similar we can prove that
$$
\alpha u_n(\alpha, z) = - s_n(\alpha,z) u_n(\alpha, z) -  \rho  s_n(\alpha, z) t_n(\alpha, z) -  z s_n(\alpha, z) + \delta_n(\alpha, z).
$$

So we have the system of equations
\begin{align}\label{eq:first}
& 1 + \alpha s_n(\alpha,z) +s_n^2(\alpha,z) = \\
& =- \frac{\rho}{2} t_n^2(\alpha,z) - \frac{z}{2} t_n(\alpha,z) - \frac{\rho}{2} u_n^2(\alpha,z) - \frac{\overline z}{2} u_n(\alpha,z) + \delta_n(\alpha,z) \nonumber \\
\label{eq:second}
& \alpha t_n(\alpha, z) = \\
&=- s_n(\alpha,z) t_n(\alpha, z) - \rho  s_n(\alpha, z) u_n(\alpha, z) - \overline z s_n(\alpha, z) + \delta_n(\alpha, z) \nonumber  \\
\label{eq:third}
& \alpha u_n(\alpha, z) = \\
&=- s_n(\alpha,z) u_n(\alpha, z) -  \rho  s_n(\alpha, z) t_n(\alpha, z) -  z s_n(\alpha, z) + \delta_n(\alpha, z) \nonumber .
\end{align}

It follows from~\eqref{eq:second} and~\eqref{eq:third} that
\begin{align*}
&(\alpha + s_n)(z t_n + \rho t_n^2) = -s_n (z \rho u_n + \overline z \rho t) -\rho^2 s_n t_n u_n - |z|^2 s_n + \delta_n(\alpha, z) \\
&(\alpha + s_n)(\overline z u_n + \rho u_n^2) = - s_n (z \rho u_n + \overline z \rho t) -\rho^2 s_n t_n u_n - |z|^2 s_n + \delta_n(\alpha, z).
\end{align*}
So, we can rewrite~\eqref{eq:first}
\begin{equation}\label{eq:new first}
1 + \alpha s_n(\alpha,z) +s_n^2(\alpha,z) + \rho^2 t_n^2(\alpha,z) + z t_n(\alpha,z) = \delta_n(\alpha,z).
\end{equation}

From equations~\eqref{eq:second} and~\eqref{eq:third} we can write equation for $t_n$
\begin{align} \label{eq:new second}
\left (\alpha + s_n - \frac{|\rho|^2 s_n^2}{\alpha + s_n} \right ) t_n = \frac{\rho z s_n^2}{\alpha+s_n} - \overline z s_n + \delta_n(\alpha, z).
\end{align}

%From~\eqref{eq:third} we can find $u_n$
%\begin{align} \label{eq:u}
%&(\alpha + s_n) u_n = -  \rho s_n u_n -  z s_n \\
%&u_n = -\frac{\rho s_n t_n +  z s_n}{\alpha + s_n} \nonumber .
%\end{align}
%Now we can put $u_n$ to~\eqref{eq:second}
%\begin{align} \label{eq:new second}
%\left (\alpha + s_n - \frac{|\rho|^2 s_n^2}{\alpha + s_n} \right ) t_n = \frac{\overline \rho z s_n^2}{\alpha+s_n)} - \overline z s_n
%\end{align}
We denote
$$
\Delta = \left (\alpha + s_n -  \frac{|\rho|^2 s_n^2}{\alpha + s_n} \right ).
$$
%We can find that
%\begin{align*}
%\overline z \rho t_n = \frac{|\rho|^2|z|^2 s_n^2}{(\alpha + s_n) \Delta} - \frac{\overline z^2 \rho s_n}{\Delta}\\
%\overline z \rho u_n = \frac{|\rho|^2|z|^2 s_n^2}{(\alpha + s_n) \Delta} - \frac{z^2 \overline \rho s_n}{\Delta}
%\end{align*}
%and
%\begin{align*}
%|\rho|^2 s_n t_n u_n = |\rho|^2 s_n \left ( \frac{\overline \rho z s_n^2}{(\alpha + s) \Delta} - \frac{\overline z s_n}{\Delta} \right ) \left (\frac{\rho \overline z s_n^2}{(\alpha + s) \Delta} - \frac{\ z s_n}{\Delta} \right )
%\end{align*}
After simple calculations we will have
\begin{align*}
& (\alpha + s_n)(z t_n + \rho t_n^2) = \\
& -s_n \left (\frac{2 \rho^2 |z|^2 s_n^2}{(\alpha + s_n) \Delta} - \frac{\overline z^2 \rho s_n}{\Delta} - \frac{z^2 \overline \rho s_n}{\Delta}\right ) \\
& -|\rho|^2 s_n \left ( \frac{\rho z s_n^2}{(\alpha + s) \Delta} - \frac{\overline z s_n}{\Delta} \right ) \left (\frac{\rho \overline z s_n^2}{(\alpha + s) \Delta} - \frac{\ z s_n}{\Delta} \right ) - |z|^2 s_n + \delta_n(\alpha, z).
\end{align*}

We denote $y_n := s_n$ and $w_n := \alpha + (\rho t_n^2 + z t_n)/y_n$. We can rewrite equations~\eqref{eq:first},~\eqref{eq:second} and~\eqref{eq:third}
\begin{align}\label{eq: main first}
& 1 + w_n y_n + y_n^2 = \delta_n(\alpha,z)\\
\label{eq: main second}
& w_n = \alpha + \frac{\rho t_n^2 + z t_n}{y_n}\\
\label{eq: main third}
& (\alpha + s_n)(z t_n + \rho t_n^2) = \\
& -s_n \left (\frac{2 \rho^2|z|^2 y_n^2}{(\alpha + y_n) \Delta} - \frac{\overline z^2 \rho y_n}{\Delta}- \frac{z^2 \rho y_n}{\Delta}\right ) - |z|^2 y_n \nonumber\\
& -|\rho|^2 y_n \left ( \frac{\rho z y_n^2}{(\alpha + y_n) \Delta} - \frac{\overline z y_n}{\Delta} \right ) \left (\frac{\rho \overline z y_n^2}{(\alpha + y_n) \Delta} - \frac{\ z y_n}{\Delta} \right )  + \delta_n(\alpha, z)\nonumber.
\end{align}

\begin{remark}
If $\rho = 0$ then we can rewrite~\eqref{eq: main first},~\eqref{eq: main second}, and~\eqref{eq: main third}
\begin{align*}
&1 + w_n y_n + y_n^2 = \delta_n(\alpha,z)\\
& w_n = \alpha + \frac{ z t_n}{y_n}\\
& (w_n - \alpha) + (w_n - \alpha)^2 y_n - |z|^2 s_n = \delta_n(\alpha,z).
\end{align*}
This equations determine the Circular law, see \cite{GotTikh2010}.
\end{remark}

We can see that the first equation~\eqref{eq: main first} doesn't depend on $\rho$.
So the first equation will be the same for all models of random matrices described in the introduction.
On the Figure~\ref{fig:semi-circular} we draw the distribution of eigenvalues of matrix $\V$ for $\rho = 0$ (Circular law case) and $\rho = 0.5$ (Elliptic law case).
%and $\rho = 1$ (Semi-circular law). We mention here that eigenvalues of $\V \J$ are singular values of $n^{-1/2}\X$ with signs $\pm$.

\begin{figure}
\begin{center}
\scalebox{.4}{\includegraphics{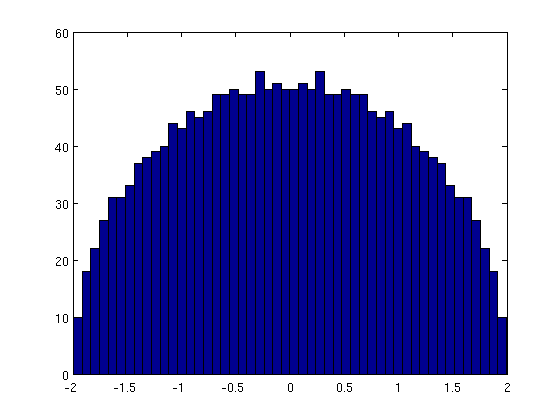}}
\scalebox{.4}{\includegraphics{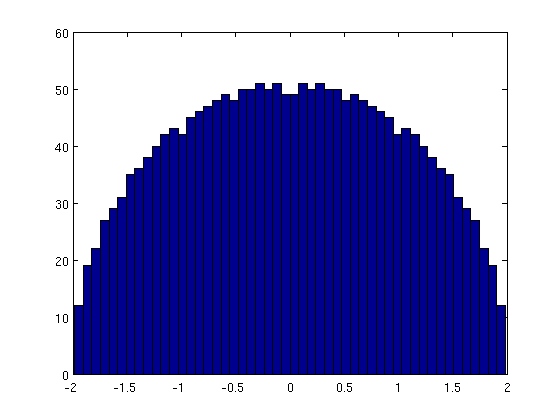}}
\end{center}
\caption{Histogram of eigenvalues of matrix $\V$ for $n = 1000$. entries are Gaussian random variables. On the left $\rho = 0$ (Circular law case). On the right $\rho = 0.5$ (Elliptic law case).}
\label{fig:semi-circular}
\end{figure}
%
%\begin{figure}[H]
%\begin{minipage}[h]{0.4\linewidth}
%\center{\includegraphics[width=1\linewidth]{circular}} a) \\
%\end{minipage}
%\begin{minipage}[h]{0.4\linewidth}
%\center{\includegraphics[width=1\linewidth]{elliptic}} b) \\
%\end{minipage}
%\begin{minipage}[h]{0.4\linewidth}
%\center{\includegraphics[width=1\linewidth]{semicir}} c) \\
%\end{minipage}
%\caption{Histogram of eigenvalues of matrix $\V \J$ for $n = 1000$, entries are Gaussian random variables with a) $\rho = 0$ (Circular law case); b) $\rho = 0.5$ (Elliptic law case); c) $\rho = 1$ (Semi-circular law case)}
%\label{fig:semi-circular}
%\end{figure}

Now we prove convergence of $s_n$ to some limit $s_0$. Let $\alpha = u + i v, v>0$. Using~\eqref{eq:new first} we write
$$
\alpha (s_n - s_m) = - (s_n - s_m)(s_n + s_m) - \rho^2(t_n - t_m)(t_n + t_m) - z (t_m - t_m) + \varepsilon_{n,m}.
$$
By triangle inequality and the fact that $|s_n| \le v^{-1}$
\begin{equation}\label{eq:sn-sm}
|s_n - s_m | \le \frac{2 |s_n - s_m|}{v^2} + \frac{\rho^2|t_n - t_m||t_n + t_m|}{v}  + \frac{|z||t_n - t_m|}{v}  + \frac{|\varepsilon_{n,m}|}{v}.
\end{equation}
From~\eqref{eq:new second} it follows that
\begin{align*}
((\alpha + s_n)^2 - \rho^2 s_n^2) t_n =  \rho z s_n^2 - \overline z \alpha s_n - \overline z s_n^2 + \varepsilon_{n}.
\end{align*}
We denote $\Delta_n : = ((\alpha + s_n)^2 - \rho^2 s_n^2)$. By triangle inequality
\begin{align}\label{eq:tn-tm}
&|\Delta_m| |t_n - t_m| \le |t_m||\Delta_n - \Delta_m| \\
&+ \frac{2|\rho||s_n - s_m| +  2 |z||s_n - s_m| }{v}  + |z||\alpha||s_n - s_m| + |\varepsilon_{n,m}| \nonumber.
\end{align}
We can find lower bound for $|\Delta_m|$:
\begin{align}\label{eq:delta}
&|\Delta_m| = |\alpha + (1-\rho) s_m||\alpha + (1 + \rho)s_m| \\
&\geq \imag(\alpha + (1-\rho) s_m)\imag(\alpha + (1+\rho) s_m) \geq v^2 \nonumber,
\end{align}
where we have used the fact that $\imag s_m \geq 0$. From definition of $\Delta_n$ it is easy to see that
\begin{equation}\label{eq:dn-dm}
|\Delta_n - \Delta_m| \le 2 |\alpha||s_n - s_m| + \frac{2 (1 + \rho^2) |s_n - s_m|}{v}.
\end{equation}
We can take $|u| \le C$, then $|\alpha| \le v + C$. From~\eqref{eq:sn-sm},\eqref{eq:tn-tm},\eqref{eq:delta} and~\eqref{eq:dn-dm} it follows that there exists constant $C'$, which depends on $\rho, C, z$,  such that
\begin{align*}
|s_n - s_m | \le \frac{C'}{v} |s_n - s_m| + |\varepsilon_{n,m}'(\alpha, z)|.
\end{align*}
We can find $v_0$ such that
$$
\frac{C'}{v} < 1 \quad \text{for all $v \geq v_0$}.
$$
Since $\varepsilon_{n,m}'(\alpha, z)$ converges to zero uniformly for all $v \geq v_0, |u| \le C$ and $s_n, s_m$ are locally bounded analytic functions in the upper half-plane we may conclude by Montel's Theorem (see \cite[Theorem~2.9]{Conw1978}) that there exists an analytic function $s_0$ in the upper half-plane such that $\lim s_n = s_0$. Since $s_n$ are
Nevanlinna functions, (that is analytic functions mapping the upper half-plane into itself) $s_0$ will be a Nevanlinna function too and there exists non-random distribution function $F(z,x)$ such that
$$
s_0(\alpha) = \int \frac{dF(z, x)}{x - \alpha}.
$$
The function $s_0$ satisfies the equations~\eqref{eq: main first},~\eqref{eq: main second}, and~\eqref{eq: main third}.
\end{proof}

\section{Acknowledgment}
The author would like to thank Prof.~Dr.~Friedrich G{\"o}tze and Prof.~Dr.~Alexander Tikhomirov for insightful comments and discussions. This research was conducted during my participation in SFB 701, Bielefeld University

\appendix
\section{Appendix}

\begin{theorem} (Central Limit Theorem) \label{th:CLT}
Let $Z_1, ... , Z_n$ be independent random variables with $\E Z_i = 0$ and finite third moment, and let $\sigma^2 = \sum_{i=1}^n \E |Z_i|^2$. Consider a standard normal variable $g$. The for every $t>0$:
$$
 \left |\mathbb P \left ( \frac{1}{\sigma} \sum_{i=1}^n Z_i \le t \right ) - \mathbb P \left ( g \le t \right ) \right | \le C \sigma^{-3} \sum_{i=1}^n \E|Z_i|^3,
$$
where $C$ is an absolute constant.
\end{theorem}

\begin{lemma}\label{l,a:decoupling}
Let event $E(X,Y)$ depends on independent random vectors $X$ and $Y$ then
$$
\Pb (E(X,Y)) \le (\Pb (E(X,Y), E(X,Y'))^{1/2},
$$
where $Y'$ is an independent copy of $Y$.
\end{lemma}
\begin{proof}
  See in~\cite{Costello2011}.
\end{proof}

\begin{lemma} \label{l,a:sum of r.v}
Let $Z_1, ... ,Z_n$ be a sequence of random variables and $p_1, ..., p_n$ be non-negative real numbers such that
$$
\sum_{i=1}^n p_i = 1,
$$
then for every $\varepsilon > 0$
$$
\Pb (\sum_{i=1}^n p_i Z_i \le \varepsilon ) \le 2 \sum_{i=1}^n p_i \Pb (Z_i \le 2 \varepsilon).
$$
\end{lemma}
\begin{proof}
  See in~\cite{Veshyn2011}.
\end{proof}

% \begin{lemma}\label{l,a:cov. comp. v}
%   Let $\mathcal{N}(T,\varepsilon)$ be an $\varepsilon$-net of set $T$. One has
%   $$
%   \mathcal{N} (Comp(\delta, \tau), 2 \tau) \le (9/\delta \tau)^{\delta n}.
%   $$
% \end{lemma}
% \begin{proof}
%   See in~\cite{Veshyn2011}.
% \end{proof}

\begin{lemma} \label{l,a:incomp vec}
If $x \in Incomp(\delta,\tau)$ then at least $\frac{1}{2}\delta \tau^2 n$ coordinates $x_k$ of $x$ satisfy
$$
\frac{\tau}{\sqrt{2 n}} \le |x_k| \le \frac{1}{\sqrt{\delta n}}.
$$
\end{lemma}
\begin{remark}\label{r:incomp vec}
We can fix some constant $c_{0}$ such that
$$
\frac{1}{4} \delta \tau^2 \le c_{0} \le \frac{1}{4}.
$$
Then for every vector $x \in Incomp(\delta, \tau)$ $|\spr(x)| = [2 c_{0} n]$.
\end{remark}
\begin{proof}
  See in~\cite{RudVesh2008}.
\end{proof}

\begin{lemma}\label{l,a:reduction}
Let $S_J = \sum_{i \in J} \xi_i$, where $J \subset [n]$, and $I \subset J$ then
$$
\sup_{v \in \R} \Pb (|S_J - v| \le \varepsilon ) \le \sup_{v \in \R} \Pb(|S_I - v| \le \varepsilon).
$$
\end{lemma}
\begin{proof}
Let us fix arbitrary $v$. From independence of $\xi_i$ we conclude
\begin{align*}
\Pb (|S_J - v| \le \varepsilon ) \le \E \Pb (|S_I + S_{J / I} - v| \le \varepsilon |\{ \xi_i\}_{i \in I}) \le \sup_{u \in \R} \Pb(|S_I - u| \le \varepsilon).
\end{align*}
\end{proof}

\begin{lemma}\label{l,a:rand var levy conc}
Let $Z$ be a random variable with $\E Z^2 \geq 1$ and with finite fourth moment, and put $M_4^4 := \E(Z-\E Z)^4$. Then for every $\varepsilon \in (0, 1)$ there exists $p = p(M_4, \varepsilon)$ such that
$$
\sup_{v \in \R} \Pb (|Z - v| \le \varepsilon) \le p.
$$
\end{lemma}
\begin{proof}
  See in~\cite{RudVesh2008}.
\end{proof}

\begin{lemma} \label{l.a:sum. c.f.}
Let $\xi_1, ... , \xi_n$ be independent random variables with $\E \xi_i^2 \geq 1$ and $\E(\xi_k - \E \xi)^4 \le M_4^4$, where $M_4$ is some finite number. Then for every $\varepsilon \in (0, 1)$ there exists $p = p(M_4, \varepsilon) \in (0, 1))$ such that the following holds: for every vector $x = (x_1, ... , x_n) \in S^{n-1}$, the sum $S = \sum_{i=1}^n x_k \xi_k$ satisfies
$$
\sup_{v \in \R} \Pb(|S-v| \le \varepsilon) \le p.
$$
\end{lemma}
\begin{proof}
  See in~\cite{RudVesh2008}.
\end{proof}

\begin{lemma}\label{l.a:tensorisation}
Let $X = (X_1, ... , X_n)$ be a random vector in $\R^n$ with independent coordinates $X_k$.\\
1. Suppose there exists numbers $\varepsilon_0 \geq 0$ and $L \geq 0$ such that
$$
\sup_{v \in \R} \Pb (|X_k - v| \le \varepsilon) \le L \varepsilon \quad \text{for all $\varepsilon \geq \varepsilon_0$ and all $k$}.
$$
Then
$$
\sup_{v \in \R^n} \Pb (||X - v||_2 \le \varepsilon \sqrt n) \le (C L \varepsilon)^n \quad \text{for all $\varepsilon \geq \varepsilon_0$},
$$
where $C$ is an absolute constant.\\
2. Suppose there exists numbers $\varepsilon > 0$ and $p \in (0, 1)$ such that
$$
\sup_{v \in \R} \Pb(|X_k - v| \le \varepsilon) \le L \varepsilon \quad \text{for all $k$}.
$$
Then there exists numbers $\varepsilon_1 = \varepsilon_1(\varepsilon, p) >0$ and $p_1 = p_1(\varepsilon, p) \in (0,1)$ such that
$$
\sup_{v \in \R^n} \Pb (||X - v||_2 \le \varepsilon_1 \sqrt n) \le (p_1)^n.
$$
\end{lemma}
\begin{proof}
  See~\cite[Lemma~3.4]{Veshyn2011}.
\end{proof}

\begin{lemma} \label{l,a:distance}
Let $1 \le m \le n$. If $\A$ has full rank, with rows $R_1, ...,R_m$ and $H = \Sp(R_j, j \neq i)$, then
$$
\sum_{i=1}^m s_i(\A)^{-2} = \sum_{i=1}^m \dist(R_i, H_i)^{-2}.
$$
\end{lemma}
\begin{proof}
  See~\cite[Lemma~A.4]{TaoVu2010}.
\end{proof}

\begin{lemma} \label{l,a:distance est}
There exist $\gamma > 0$ and $\delta > 0$ such that for all $n \gg 1$ and $1 \le i \le n$, any deterministic vector $v \in \mathbb C$ and any subspace $H$ of $\mathbb C^n$ with $1 \le \dim(H) \le n - n^{1-\gamma}$, we have, denoting $R:=(X_1,...,X_n) + v$,
$$
\mathbb P( \dist(R,H) \le \frac{1}{2} \sqrt{n - \dim(H)}) \le \exp(-n^\delta).
$$
\end{lemma}
\begin{proof}
    See~\cite[Statement~5.1]{TaoVu2010}.
\end{proof}

\begin{lemma} \label{stilt.transform}
Under the condition $\Cond$ for $\alpha = u + i v, v > 0$
$$
\E \left | \frac{1}{n} \sum_{i=1}^n R_{i i}(\alpha, z) - \E \left (\frac{1}{n} \sum_{i=1}^n R_{ii} (\alpha,z) \right ) \right |^2 \le \frac{C}{n v^2}.
$$
\end{lemma}
\begin{proof}
To prove this lemma we will use Girko's method. Let $\X^{(j)}$ be a matrix $\X$ with $j$-th row and column removed. Define matrices $\V^{(j)}$ and $\V^{(j)}(z)$
as in~\eqref{eq:matrices} and $\RR^{(j)}$ by~\eqref{eq:resolvent}. It is easy to see that
$$
\Rank(\V(z) - \V^{(j)}(z)) = \Rank ( \V- \V^{(j)}) \le 4.
$$
Then
\begin{equation} \label{eq:rank}
\frac{1}{n}| \Tr (\V(z)  - \alpha \I)^{-1} - \Tr (\V^{(j)}(z) - \alpha \I)^{-1}| \le \frac{\Rank (\V(z) - \V^{(j)}(z))}{n v} \le \frac{4}{n v}.
\end{equation}

We introduce the family of $\sigma$-algebras $\mathcal{F}_i = \sigma \{ X_{j,k}, j, k > i \}$ and conditional mathematical expectation $\E_i = \E ( \cdot |\mathcal{F}_i)$
with respect to this $\sigma$-algebras. We can write
\begin{align*}
\frac{1}{n} \Tr \RR - \frac{1}{n} \E \Tr \RR  = \frac{1}{n} \sum_{i = 1}^n \E_i \Tr \RR - \E_{i-1} \Tr \RR = \sum_{i = 1}^n \gamma_i.
\end{align*}
The sequence $(\gamma_i, \mathcal{F}_i)_{i \geq 1}$ is a martingale difference. By~\eqref{eq:rank}
\begin{align} \label{eq:gamma}
& |\gamma_i| =  \frac{1}{n} | \E_i (\Tr \RR - \Tr \RR^{(i)}) - \E_{i-1} (\Tr \RR - \Tr \RR^{(i)})| \le \\
& \le | \E_i (\Tr \RR - \Tr \RR^{(i)})|  + |\E_{i-1} (\Tr \RR - \Tr \RR^{(i)})| \le \frac{C}{v n}.
\end{align}
From Burkholder inequality for martingale difference (see~\cite{Shiryaev1996})
$$
\E \left |\sum_{i=1}^n \gamma_i \right |^2 \le K_2 \E \left (\sum_{i=1}^n |\gamma_i|^2 \right )
$$
and~\eqref{eq:gamma} it follows
$$
\E \left | \frac{1}{n} \sum_{i=1}^n R_{i i}(\alpha, z) - \E \left (\frac{1}{n} \sum_{i=1}^n R_{ii} (\alpha,z) \right ) \right |^2 \le K_2 \E (\sum_{i=1}^n |\gamma_i|^2) \le K_2 \frac{C}{n v^2}.
$$
\end{proof}

\begin{lemma}\label{l,a:second diag}
Under the condition $\Cond$ for $\alpha = u + i v, v > 0$
$$
\E \left | \frac{1}{n} \sum_{i=1}^n R_{i,i+n}(\alpha, z) - \E \left (\frac{1}{n} \sum_{i=1}^n R_{i,i+n} (\alpha,z) \right ) \right |^2 \le \frac{C}{n v^4}.
$$
\end{lemma}

\begin{proof}
As in Lemma~\ref{stilt.transform} we introduce matrices $\V^{(j)}$ and $\RR^{(j)}$. We have
$$
\V = \V^{(j)} + e_j e_j^T \V + \V e_j e_j^T + e_{j+n} e_{j+n}^T \V  + \V e_{j+n} e_{j+n}^T
$$
By resolvent equality $\RR - \RR^{(j)} = -\RR^{(j)} (\V(z) - \V^{(j)}(z)) \RR$
\begin{align*}
&\frac{1}{n}\sum_{k = 1}^n (\RR_{k, k+n} - \RR_{k, k+n}^{(j)}) = \\
& = \frac{1}{n} \sum_{k=1}^n [\RR^{(j)} (e_j e_j^T \V + e_{j+n}e_{j+n}^T \V + \V e_j e_j^T  + \V  e_{j+n}e_{j+n}^T) \RR ]_{k,k+n} =\\
& = \mathbb T_1 + \mathbb T_2 + \mathbb T_3 + \mathbb T_4.
\end{align*}
Let us consider the first term. The arguments for other terms are similar.
$$
\sum_{k=1}^n [\RR^{(j)} e_j e_j^T \V \RR]_{k,k+n} = \Tr \RR^{(j)} e_j e_{j}^T \V \RR \EE =
\sum_{i=1}^{2n} [\RR \EE \RR^{(j)}]_{ij}[e_j e_{j}^T \V]_{j i},
$$
where
$$
\EE = \begin{pmatrix}
\OO_n & \OO_n\\
\I &  \OO_n
\end{pmatrix}.
$$
From $\max (||\RR^{(j)}||, ||\RR||) \le v^{-1}$ and H\"older inequality it follows that
$$
\E \left |\sum_{k=1}^n [\RR^{(j)} e_j e_j^T \V \J \RR]_{k,k+n} \right |^2 \le \frac{C}{v^4}.
$$
By similar arguments as in Lemma~\ref{stilt.transform} we can conclude the statement of the Lemma.
\end{proof}

\begin{lemma}\label{l,a:stein}
 Under the condition $\Cond$ for $\alpha = u + i v, v > 0$
\begin{align*}
&\frac{1}{n^{3/2}}\sum_{j,k=1}^n \E X_{j k} R_{k+n,j} = \\
& = \frac{1}{n^{2}}\sum_{j,k=1}^n \E \left [ \frac{ \partial \RR}{\partial X_{j k}} \right ]_{k+n,j} +
 \frac{\rho}{n^{2}}\sum_{j,k=1}^n \E \left [ \frac{ \partial \RR}{\partial X_{k j}} \right ]_{k+n,j}  + r_n(\alpha,z),
\end{align*}
where
$$
|r_n(\alpha,z)| \le \frac{C}{\sqrt n v^3}
$$
\end{lemma}

\begin{proof}
By Taylor's formula
\begin{align} \label{eq:taylor 1}
&\E X f(X,Y) =  f(0,0) \E X +  f_x'(0,0) \E X^2 + f_y'(0,0) \E X Y +\\
&+\E (1 - \theta) [X^3 f_{x x}''(\theta X, \theta Y) + 2 X^2 Y f_{x y}''(\theta X, \theta Y) + X Y^2 f_{y y}''(\theta X, \theta Y) ] \nonumber
\end{align}
and
\begin{align} \label{eq:taylor 2}
&\E f_x'(X,Y) = f_x'(0,0) + \E (1-\theta) [X f_{x x}''(\theta X, \theta Y) + Y f_{x y}''(\theta X, \theta Y)]\\
&\E f_y'(X,Y) = f_y'(0,0) + \E (1-\theta) [X f_{x y}''(\theta X, \theta Y) + Y f_{y y}''(\theta X, \theta Y)], \nonumber
\end{align}
where $\theta$ has uniform distribution on $[0, 1]$. From~\eqref{eq:taylor 1} and~\eqref{eq:taylor 2} for $j \neq k$
\begin{align*}
&\left |\E X_{j k} R_{k+n,j} - \E \left [ \frac{ \partial \RR}{\partial X_{j k}} \right ]_{k+n,j} -
 \rho \E \left [ \frac{ \partial \RR}{\partial X_{k j}} \right ]_{k+n,j} \right | \le \\
& (|X_{j k}|^3 + |X_{j k}|)  \left | \left [\frac{ \partial^2 \RR}{\partial X_{j k}^2} (\theta X_{j k}, \theta X_{k j}) \right ]_{k+n, j} \right | + \\
&(|X_{k j}|^2|X_{j k}| + |X_{k j}|) \left | \left [ \frac{ \partial^2 \RR}{\partial X_{k j}^2} (\theta X_{j k}, \theta X_{k j}) \right ]_{k+n, j} \right | + \\
& (2|X_{j k}|^2 |X_{k j}| + |X_{j k}| + |X_{k j}| ) \left | \left [ \frac{ \partial^2 \RR}{\partial X_{j k} \partial X_{k j}} (\theta X_{j k}, \theta X_{k j}) \right ]_{k+n, j} \right | .
\end{align*}
Let us consider the first term in the sum. The bounds for the second and third terms can be obtained by similar arguments. We have
\begin{align*}
& \frac{\partial^2 \RR}{\partial X_{j k}^2} = \frac{1}{n} \RR (e_j e_{n+k}^T + e_{n+k}e_j^T) \RR (e_j e_{n+k}^T + e_{n+k}e_j^T) \RR
& = \mathbb P_1 + \mathbb P_2 + \mathbb P_3 + \mathbb P_4,
\end{align*}
where
\begin{align*}
&\mathbb P_1 = \frac{1}{n} \RR e_j e_{n+k}^T \RR e_j e_{n+k}^T \RR \\
&\mathbb P_2 = \frac{1}{n} \RR e_j e_{n+k}^T \RR e_{n+k} e_j^T \RR \\
&\mathbb P_3 = \frac{1}{n} \RR e_{n+k} e_j^T \RR e_j e_{n+k}^T \RR  \\
&\mathbb P_4 = \frac{1}{n} \RR e_{n+k} e_j^T \RR e_{n+k} e_j^T \RR.
\end{align*}
From $|\RR_{i,j}| \le v^{-1}$ it follows that
$$
\frac{1}{n^{5/2}} \sum_{j,k =1}^n \E |X_{j k}|^{\alpha} |[\mathbb P_i]_{n+k,j}| \le \frac{C}{\sqrt n v^3}
$$
for $\alpha = 1, 3$ and $ i = 1, ..., 4$. For $j = k$
$$
\frac{1}{n^{2}} \sum_{j = 1}^n \E \left [\frac{ \partial \RR}{\partial X_{j j}} \right ]_{j+n,j} = \frac{1}{n^{2}} \sum_{j = 1}^ n ( \E R_{j+n, j}^2 + \E R_{j, j}R_{j+n, j+n} ) \le \frac{C}{n v^2}.
$$
So we can add this term to the sum
$$
\frac{\rho}{n^{2}}\sum_{\substack {j,k=1\\ j \neq k} }^n \E \left [ \frac{ \partial \RR}{\partial X_{k j}} \right ]_{k+n,j}.
$$
\end{proof}
\newpage
%\section{Figures}
%\begin{figure}[H]
%\begin{minipage}[h]{1\linewidth}
%\center{\includegraphics[width=0.4\linewidth]{1g}} a) \\
%\end{minipage}
%\begin{minipage}[h]{1\linewidth}
%\center{\includegraphics[width=0.4\linewidth]{1b}} b) \\
%\end{minipage}
%
%\begin{minipage}[h]{1\linewidth}
%\center{\includegraphics[width=0.4\linewidth]{3g}} c) \\
%\end{minipage}
%\begin{minipage}[h]{1\linewidth}
%\center{\includegraphics[width=0.4\linewidth]{2g}} d) \\
%\end{minipage}
%
%\caption{Spectrum of matrix $n^{-1/2} \X$ for $n = 2000$,
%a) entries are Gaussian random variables with $\rho = 0.5$
%b) entries are Rademacher random variables with $\rho = 0.5$;
%c) entries are Gaussian random variables with $\rho = -0.5$
%d) entries are Gaussian random variables with $\rho = 0$  }
%\label{fig:ellip and circ}
%\end{figure}
%
%
%\begin{figure}[H]
%\begin{minipage}[h]{0.4\linewidth}
%\center{\includegraphics[width=1\linewidth]{circular}} a) \\
%\end{minipage}
%\begin{minipage}[h]{0.4\linewidth}
%\center{\includegraphics[width=1\linewidth]{elliptic}} b) \\
%\end{minipage}
%\begin{minipage}[h]{0.4\linewidth}
%\center{\includegraphics[width=1\linewidth]{semicir}} c) \\
%\end{minipage}
%\caption{Histogram of eigenvalues of matrix $\V \J$ for $n = 1000$, entries are Gaussian random variables with a) $\rho = 0$ (Circular law case); b) $\rho = 0.5$ (Elliptic law case); c) $\rho = 1$ (Semi-circular law case)}
%\label{fig:semi-circular}
%\end{figure}
%\newpage

%
\bibliographystyle{plain}
\bibliography{literatur}

\end{document}